\documentclass[A4j,11pt]{article}
\usepackage{graphics}
\usepackage{amsfonts,amsmath,amscd}
\usepackage{amssymb}
\begin{document}

\title{{\bf Homogeneity of infinite dimensional \\
anti-Kaehler isoparametric submanifolds}}
\author{{\bf Naoyuki Koike}}

\date{}

\maketitle

{\small\textit{Department of Mathematics, Faculty of Science, 
Tokyo University of Science,}}

{\small\textit{1-3 Kagurazaka Shinjuku-ku, Tokyo 162-8601 Japan}}

{\small\textit{E-mail address}: koike@ma.kagu.tus.ac.jp}

%
%
%

\begin{abstract}
In this paper, we prove that, if a full irreducible infinite dimensional anti-Kaehler 
isoparametric submanifold of codimension greater than one has $J$-diagonalizable 
shape operators, then it is homogeneous.  
\end{abstract}

\vspace{0.1truecm}

{\small\textit{Keywords}: $\,$ anti-Kaehler isoparametric submanifold, $J$-principal 
curvature,

\hspace{1.9truecm}J-curvature distribution, regularizability}




\vspace{0.3truecm}

\section{Introduction}
In 1999, E. Heintze and X. Liu [HL2] proved that all irreducible 
isoparametric submanifolds of codimension greater than one in the (separable) 
Hilbert space are homogeneous, which is the infinite 
dimensional version of the homogeneity theorem for isoparametric submanifolds 
in a (finite dimensional) Euclidean space by G. Thorbergsson ([Th]).  Note that 
the result of Thorbergsson states that all irreducible isoparametric 
submanifolds of codimension greater than two in a Euclidean space are 
homogeneous.  In 2002, by using this result of Heintze-Liu, 
U. Christ [Ch] proved that all irreducible equifocal submanifolds with flat 
section of codimension greater than one in a simply connected symmetric space of compact 
type are homogeneous, where we note that, in a simply connected symmetric space of compact 
type, the notion of an equifocal submanifold coincides with that of an isoparametric 
submanifold with flat section in the sense of [HLO].  
In [Koi1], we introduced the notion of a complex equifocal submanifold in a symmetric space 
of non-compact type.  Here we note that all isoparametric submanifolds with flat section 
are complex equifocal.  In [Koi2], we showed that the study of complex equifocal 
$C^{\omega}$-submanifolds in the symmetric spaces are reduced to that of anti-Kaehler 
isoparametric submanifolds in the infinite dimensional anti-Kaehler space, where 
$C^{\omega}$ means the real analyticity.  
In this paper, we shall investigate an anti-Kaehler isoparametric 
submanifold with $J$-diagonalizable shape opeartors.  According to the discussion in [Koi2], 
we can show that the study of certain kind of isoparametric submanifolds 
with flat section in symmetric spaces of non-compact type are reduced to that of 
anti-Kaehler isoparametric submanifolds with $J$-diagonalizable shape operators 
in the infinite dimensional anti-Kaehler space, which is called a proper anti-Kaehler 
isoparametric submanifold in [Koi2].  
L. Geatti and C. Gorodski ([GG]) introduced the notion of an isoparametric submanifold with 
diagonalizable Weingarten operators in a finite dimensional pseudo-Euclidean space.  
Note that anti-Kaehler isoparametric submanifolds with 
$J$-diagonalizable shape operators give a subclass of the class of the infinite dimensional 
version of isoparametric submanifolds with diagonalizable Weingarten operators 
(see Remark 2.1).  

In this paper, we prove the following homogeneity theorem for 
anti-Kaehler isoparametric $C^{\omega}$-submanifolds with $J$-diagonalizable 
shape operators in the infinite dimensional anti-Kaehler space.  

\vspace{0.5truecm}

\noindent
{\bf Theorem A.} {\sl Let $M$ be a full irreducible anti-Kaehler isoparametric 
$C^{\omega}$-submanifold with $J$-diagonalizable shape operators of codimension greater than 
one in the infinite dimensional anti-Kaehler space.  Then $M$ is homogeneous.}

\vspace{0.5truecm}

\noindent
{\it Remark 1.1.} This homogeneity theorem will be useful to prove 
homogeneity of certain kind of isoparametric submanifolds with flat 
section in symmetric spaces of non-compact type, which have principal orbits of Hermann 
actions as homogeneous examples.  

\section{Basic notions and facts}
In this section, we shall first recall the notion of an anti-Kaehler isoparametric 
submanifold in the infinite dimensional anti-Kaehler space introduced in [Koi2].  
Let $V$ be an infinite dimensional topological real vector space (or a finite dimensional 
real vector space), $\widetilde J$ a continuous linear operator of $V$ such that 
$\widetilde J^2=-{{\rm id}}$ and $\langle\,\,,\,\,\rangle$ a continuous non-degenerate 
symmetric bilinear form of $V$ 
such that $\langle\widetilde JX,\widetilde JY\rangle=-\langle X,Y\rangle$ 
holds for every $X,Y\in V$.  Assume that there exists an orthogonal time-space 
decomposition $V=V_-\oplus V_+$ (i.e., $\langle\,\,,\,\,\rangle\vert_{V_-\times V_+}=0$, 
$\langle\,\,,\,\,\rangle\vert_{V_-\times V_-}:$ negative definite, 
$\langle\,\,,\,\,\rangle\vert_{V_+\times V_+}:$ positive definite) 
such that $\widetilde JV_-=V_+,\,\,(V,\langle\,\,,\,\,\rangle_{V_{\pm}})$ is a separable 
Hilbert space and that the distance topology associated with 
$\langle\,\,,\,\,\rangle_{V_{\pm}}$ coincides with the original topology of $V$, where 
$\langle\,\,,\,\,\rangle_{V_{\pm}}:=
-\pi_{V_-}^{\ast}\langle\,\,,\,\,\rangle+\pi_{V_+}^{\ast}\langle\,\,,\,\,
\rangle$ ($\pi_{V_{\pm}}\,:\,$ the projection of $V$ onto $V_{\pm}$).  
Then we call $(V,\langle\,\,,\,\,\rangle,\widetilde J)$ the {\it anti-Keahler space}.  
Let $M$ be a Hilbert manifold modelled on a separable Hilbert space 
$(V',\langle\,\,,\,\,\rangle_{V'})$.  
Let $\langle\,\,,\,\,\rangle$ be a section of the $(0,2)$-tensor bundle 
$T^{\ast}M\otimes T^{\ast}M$ such that $\langle\,\,,\,\,\rangle_x$ is 
a continuous non-degenerate symmetric bilinear form on $T_xM$ for each 
$x\in M$ and $J$ a section of the $(1,1)$-tensor bundle 
$T^{\ast}M\otimes TM$ such that 
$J^2=-{\rm id},\,\,\nabla J=0$ ($\nabla\,:\,$ the Levi-Civita connection of 
$\langle\,\,,\,\,\rangle$), $J_x$ is a continuous linear operator of 
$T_xM$ for each $x\in M$ and that $\langle JX,JY\rangle=-\langle X,Y\rangle$ 
for every $X,Y\in TM$.  We call $(M,\langle\,\,,\,\,\rangle,J)$ an 
{\it anti-Keahler Hilbert manifold} if, for each $x\in M$, there exist 
distributions $W_{\pm}$ on some neighborhood $U$ of $x$ satisfying 
the following condition:

\vspace{0.2truecm}

For each $y\in U$, $(W_{\pm})_y$ gives an orthogonal time-space 
decomposition of $(T_yM,\langle\,\,,\,\,\rangle_y)$ 

(i.e., $T_yM=(W_-)_y\oplus(W_+)_y,\,\langle\,\,,\,\,\rangle_y\vert_{(W_-)_y\times(W_+)_y}=0, 
\langle\,\,,\,\,\rangle_y\vert_{(W_-)_y\times(W_-)_y}\,:\,$negative 

definite and $\langle\,\,,\,\,\rangle_y\vert_{(W_+)_y\times(W_+)_y}\,:\,$positive definite), 
$(T_yM,\langle\,\,,\,\,\rangle_{y,(W_{\pm})_y})$ is isometric 

to $(V',\langle\,\,,\,\,\rangle_{V'})$ and $J_y(W_-)_y=(W_+)_y$, where 
$\langle\,\,,\,\,\rangle_{y,(W_{\pm})_y}:=
-\pi_{(W_-)_y}^{\ast}\langle\,\,,\,\,\rangle_y+$

$\pi_{(W_+)_y}^{\ast}\langle\,\,,\,\,\rangle_y$ ($\pi_{(W_{\pm})_y}\,:\,$ the projection of 
$T_yM$ onto $(W_{\pm})_y$).  

\vspace{0.2truecm}

\noindent
Let $f$ be an isometric immersion of an anti-Keahler Hilbert manifold 
$(M,\langle\,\,,\,\,\rangle_M,J)$ into an anti-Keahler space 
$(V,\langle\,\,,\,\,\rangle,\widetilde J)$.  If $f^{\ast}\langle\,\,,\,\,\rangle
=\langle\,\,,\,\,\rangle_M$ and if $\widetilde J\circ f_{\ast}=f_{\ast}\circ J$ holds, 
then we call $(M,\langle\,\,,\,\,\rangle_M,J)$ (or $M$) an {\it anti-Kaehler submanifold in} 
$(V,\langle\,\,,\,\,\rangle,\widetilde J)$ {\it immersed by} $f$.  
If $M$ is of finite codimension and, 
for each $v\in T^{\perp}M$, the shape operator $A_v$ is a compact operator 
with respect to $f^{\ast}\langle\,\,,\,\,\rangle_{V_{\pm}}$, then we call 
$(M,\langle\,\,,\,\,\rangle_M,J)$ (or $M$) an {\it anti-Kaehler Fredholm submanifold}.  
Let $M$ be an anti-Kaehler Fredholm submanifold.  Denote by $A$ the shape tensor of $M$.  
Fix a unit normal vector $v$ of $M$.  If there exists $X(\not=0)\in TM$ with 
$A_vX=aX+bJX$, then we call the complex number $a+b\sqrt{-1}$ a 
$J$-{\it eigenvalue of} $A_v$ 
(or a $J$-{\it principal curvature of direction} $v$) 
and call $X$ a $J$-{\it eigenvector for} $a+b\sqrt{-1}$.  
Also, we call the space of all $J$-eigenvectors for 
$a+b\sqrt{-1}$ a $J$-{\it eigenspace for} $a+b\sqrt{-1}$.  
The $J$-eigenspaces are orthogonal to one another and 
they are $J$-invariant, respectively.  
We call the set of all $J$-eigenvalues of $A_v$ the $J$-{\it spectrum of} 
$A_v$ and denote it by ${{\rm Spec}}_JA_v$.  
Since $M$ is an anti-Kaehler Fredholm submanifold, the set 
${{\rm Spec}}_JA_v\setminus\{0\}$ is described as follows:
$${{\rm Spec}}_JA_v\setminus\{0\}
=\{\mu_i\,\vert\,i=1,2,\cdots\}$$
$$\left(
\begin{array}{c}
\displaystyle{\vert\mu_i\vert>\vert\mu_{i+1}\vert\,\,\,{{\rm or}}
\,\,\,{\rm "}\vert\mu_i\vert=\vert\mu_{i+1}\vert\,\,\&\,\,
{{\rm Re}}\,\mu_i>{{\rm Re}}\,\mu_{i+1}{\rm "}}\\
\displaystyle{{{\rm or}}\,\,\,"\vert\mu_i\vert=\vert\mu_{i+1}\vert 
\,\,\&\,\,
{{\rm Re}}\,\mu_i={{\rm Re}}\,\mu_{i+1}\,\,\&\,\,
{{\rm Im}}\,\mu_i=-{{\rm Im}}\,\mu_{i+1}>0"}
\end{array}
\right).$$
Also, the $J$-eigenspace for each $J$-eigenvalue of $A_v$ other than $0$ 
is of finite dimension.  We call the $J$-eigenvalue $\mu_i$ 
the $i$-{\it th} $J$-{\it principal curvature of direction} $v$.  
Assume that the normal holonomy group of $M$ is trivial.  
Fix a parallel normal vector field $\widetilde v$ of $M$.  
Assume that the number (which may be $\infty$) of distinct $J$-principal 
curvatures of direction $\widetilde v_x$ is independent of the choice of $x\in M$.  
Then we can define complex-valued functions $\widetilde{\mu}_i$ 
($i=1,2,\cdots$) on $M$ by assigning the $i$-th $J$-principal curvature of direction 
$\widetilde v_x$ to each $x\in M$.  We call this function $\widetilde{\mu}_i$ the 
$i$-{\it th} $J$-{\it principal curvature function of direction} $\widetilde v$.  
The submanifold $M$ is called an {\it anti-Kaehler isoparametric submanifold} if it 
satisfies the following condition:

\vspace{0.2truecm}

The normal holonomy group of $M$ is trivial, and, for each parallel normal vector field 

$\widetilde v$ of $M$, the number of distinct $J$-principal curvatures of direction 
$\widetilde v_x$ is independent 

of the choice of $x\in M$, each $J$-principal curvature function of direction $\widetilde v$ 
is constant 

on $M$ and it has constant multiplicity.  

\vspace{0.2truecm}

\noindent
Let $M$ be an anti-Kaehler Fredholm submanifold in $V$.  
Let $\{e_i\}_{i=1}^{\infty}$ be an orthonormal system of $T_xM$.  If 
$\{e_i\}_{i=1}^{\infty}\cup\{Je_i\}_{i=1}^{\infty}$ is an orthonormal base 
of $T_xM$, then we call $\{e_i\}_{i=1}^{\infty}$ (rather than 
$\{e_i\}_{i=1}^{\infty}\cup\{Je_i\}_{i=1}^{\infty}$) a 
$J$-{\it orthonormal base}.  
If there exists a $J$-orthonormal base consisting of $J$-eigenvectors of 
$A_v$, then we say that $A_v$ {\it is diagonalized with respect to a $J$-orthonormal base} 
(or $A_v$ {\it is} $J$-{\it diagonalizable}).  
If, for each $v\in T^{\perp}M$, the shape operator $A_v$ is $J$-diagonalizable, then 
we say that $M$ {\it has} $J$-{\it diagonalizable shape operators}.  

\vspace{0.5truecm}

\noindent
{\it Remark 2.1.} 
If $A_v$ is diagonalized with respect to a $J$-orthonormal base, then 
the complexification $A_v^{\bf c}$ of $A_v$ is diagonalized with respect to an orthonormal 
base.  In fact, if $A_vX=aX+bJX$, then we have 
$A_v^{\bf c}(X\pm\sqrt{-1}JX)=(a\mp\sqrt{-1}b)(X\pm\sqrt{-1}JX)$.  

\vspace{0.5truecm}

\noindent
Let $M$ be an anti-Kaehler isoparametric submanifold with $J$-diagonalizable shape 
operators, where we note that such a submanifold was called 
a {\it proper anti-Kaehler isoparametric submanifold} in [Koi2] 
(in this paper, we do not use this terminology).  
Then, since the ambient space is flat and the normal holonomy group of $M$ is 
trivial, it follows from the Ricci equation that the shape operators $A_{v_1}$ 
and $A_{v_2}$ commute for arbitrary two normal vector $v_1$ and $v_2$ 
of $M$.  Hence the shape operators $A_v$'s ($v\in T^{\perp}_xM$) 
are simultaneously diagonalized with respect to a $J$-orthonormal base.  
Let $\{E_0\}\cup\{E_i\,\vert\,i\in I\}$ be the family of distributions on $M$ such that, 
for each $x\in M$, 
$\{(E_0)_x\}\cup\{(E_i)_x\,\vert\,i\in I\}$ is the set of all common $J$-eigenspaces of 
$A_v$'s ($v\in T^{\perp}_xM$), where 
$\displaystyle{(E_0)_x=\mathop{\cap}_{v\in T^{\perp}_xM}{\rm Ker}\,A_v}$.  
For each $x\in M$, $T_xM$ is equal to the closure 
$\overline{\displaystyle{(E_0)_x\oplus
\left(\mathop{\oplus}_{i\in I}(E_i)_x\right)}}$ of 
$\displaystyle{(E_0)_x\oplus\left(\mathop{\oplus}_{i\in I}(E_i)_x\right)}$.  
We regard $T^{\perp}_xM$ ($x\in M$) as a complex vector space by 
$J_x\vert_{T^{\perp}_xM}$ and denote the dual space of the complex vector 
space $T^{\perp}_xM$ by $(T^{\perp}_xM)^{\ast_{\bf c}}$.  
Also, denote by $(T^{\perp}M)^{\ast_{\bf c}}$ the complex vector bundle over 
$M$ having $(T^{\perp}_xM)^{\ast_{\bf c}}$ as the fibre over $x$.  
Let $\lambda_i$ ($i\in I$) be the section of $(T^{\perp}M)^{\ast_{\bf c}}$ 
such that $A_v={{\rm Re}}(\lambda_i)_x(v){{\rm id}}+{{\rm Im}}(\lambda_i)_x(v)
J_x$ on $(E_i)_x$ for any $x\in M$ and any $v\in T^{\perp}_xM$.  
We call $\lambda_i$ ($i\in I$) $J$-{\it principal curvatures} of $M$ and 
$E_i$ ($i\in I$) $J$-{\it curvature distributions} of $M$.  
The distribution $E_i$ is integrable and each leaf of $E_i$ is a complex 
sphere.  Each leaf of $E_i$ is called a {\it complex curvature sphere}.  
It is shown that there uniquely exists a normal vector field $n_i$ of $M$ with 
$\lambda_i(\cdot)
=\langle n_i,\cdot\rangle-\sqrt{-1}\langle Jn_i,\cdot\rangle$ 
(see Lemma 5 of [Koi2]).  We call $n_i$ ($i\in I$) the 
$J$-{\it curvature normals} of $M$.  Note that $n_i$ is parallel with 
respect to the normal connection of $M$.  
Set ${\it l}^x_i:=(\lambda_i)_x^{-1}(1)$.  According to (i) of Theorem 2 in 
[Koi2], the tangential focal set of $M$ at $x$ is equal to 
$\displaystyle{\mathop{\cup}_{i\in I}{\it l}_i^x}$.  
We call each ${\it l}_i^x$ a {\it complex focal hyperplane of} $M$ {\it at} 
$x$.  Let $\widetilde v$ be a parallel normal vector field of $M$.  
If $\widetilde v_x$ belongs to at least one ${\it l}_i$, then it is called 
a {\it focal normal vector field} of $M$.  
For a focal normal vetor field $\widetilde v$, the focal map 
$f_{\widetilde v}$ is defined by $f_{\widetilde v}(x):=x+\widetilde v_x\,\,\,(x\in M)$.  The image $f_{\widetilde v}(M)$ is called a {\it focal submanifold} of $M$, which we denote by 
$F_{\widetilde v}$.  
For each $x\in F_{\widetilde v}$, the inverse image $f_{\widetilde v}^{-1}(x)$ 
is called a {focal leaf} of $M$.  
Denote by $T_i^x$ the complex reflection of order $2$ with respect to 
${\it l}_i^x$ (i.e., the rotation of angle $\pi$ having ${\it l}_i^x$ as the 
axis), which is an affine transformation of $T^{\perp}_xM$.  
Let ${\cal W}_x$ be the group generated by $T_i^x$'s ($i\in I$).  
According to Proposition 3.7 of [Koi3], ${\cal W}_x$ is discrete.  
Furthermore, it follows from this fact that ${\cal W}_x$ is isomorphic to 
an affine Weyl group.  This group ${\cal W}_x$ is independent of the choice 
of $x\in M$ (up to group isomorphicness).  
Hence we simply denote it by ${\cal W}$.  We call this group the 
{\it complex Coxeter group associated with} $M$.  
According to Lemma 3.8 of [Koi3], 
$W$ is decomposable (i.e., it is decomposed into a non-trivial product of 
two discrete complex reflection groups) if and only if there exist two 
$J$-invariant linear subspaces $P_1$ ($\not=\{0\}$) and $P_2$ ($\not=\{0\}$) 
of $T^{\perp}_xM$ such that $T^{\perp}_xM=P_1\oplus P_2$ (orthogonal 
direct sum), $P_1\cup P_2$ contains all $J$-curvature normals of 
$M$ at $x$ and that $P_i$ ($i=1,2$) contains at least one $J$-
curvature normal of $M$ at $x$.  
Also, according to Theorem 1 of [Koi3], $M$ is irreducible if and only if 
${\cal W}$ is not decomposable.  

Next we shall recall the notion of an aks-representation.  
Let $(N,\langle\,\,,\,\,\rangle,J)$ be a finite dimensional anti-Keahler manifold.  
If there exists an involutive holomorphic isometry $s_p$ of $N$ having $p$ as an isolated 
fixed point for each $p\in N$, then we call $(N,J,\langle\,\,,\,\,\rangle)$ an 
{\it anti-Keahler symmetric space}.  Furthermore, if the isometry group of 
$(N,J,\langle\,\,,\,\,\rangle)$ is semi-simple, then it is said to be semi-simple.  Let 
$G$ be a connected complex Lie group and $K$ a closed complex subgroup of $G$.  
If there exists an involutive complex automorphism $\rho$ of $G$ such that 
$G_{\rho}^0\subset K\subset G_{\rho}$ ($G_{\rho}\,:\,$ the group 
of all fixed points of $\rho$, $G_{\rho}^0\,:\,$ the identity component of $G_{\rho}$) 
, then we call the pair $(G,K)$ an {\it anti-Keahler symmetric pair}.  
We [Koi4] showed that, for each anti-Kaehler symmetric pair $(G,K)$, the quotient $G/K$ 
is an anti-Kaehler symmetric space in a natural manner and that, conversely, from each 
anti-Kaehler symmetric space, an anti-Kaehler symmetric pair arises.  
Let $\mathfrak g$ be a complex Lie algebra and $\tau$ a complex involution 
of $\mathfrak g$.  Then we call $(\mathfrak g,\tau)$ an 
{\it anti-Kaehler symmetric Lie algebra}.  
We [Koi4] showed that an anti-Kaehler symmetric Lie algebra arises from an 
anti-Kaehler symmetric pair and that, conversely, an anti-Kaehler symmetric pair 
arises from an anti-Kaehler symmetric Lie algebra.  
Let $(N,J,\langle\,\,,\,\,\rangle)$ be an irreducible 
anti-Kaehler symmetric space, $G$ the identity component of the holomorphic isometry 
group of $(N,J,\langle\,\,,\,\,\rangle)$ and $K$ the isotropy group of $G$ at 
some point $x_0\in N$, where the irreducibility implies that $N$ is 
not decomposed into the non-trivial product of two anti-Kaehler symmetric 
spaces.  Assume that $(N,J,\langle\,\,,\,\,\rangle)$ 
does not have the pseudo-Euclidean part in its de Rham decomposition.  
Note that an anti-Kaehler symmetric space without pseudo-Euclidean part 
is not necessarily semi-simple (see [CP],[W1]).  
Let $G/K$ be an irreducible anti-Kaehler symmetric space and $(\mathfrak g,\tau)$ the 
anti-Kaehler symmetric Lie algebra associated with $G/K$.  
Also, set $\mathfrak p:={\rm Ker}(\tau+{\rm id})$.  
The space ${\rm Ker}(\tau-{\rm id})$ is equal to the Lie algebra $\mathfrak k$ of $K$ 
and $\mathfrak p$ is identified with $T_{eK}(G/K)$.  
Denote by ${\rm Ad}_G$ be the adjoint representation of $G$.  Define 
${\rm Ad}_G\vert_{\mathfrak p}:K\to{\rm GL}(\mathfrak p)$ by 
$({\rm Ad}_G\vert_{\mathfrak p})(k):={\rm Ad}_G(k)\vert_{\mathfrak p}\,\,(k\in K)$.  
We call this representation ${\rm Ad}_G\vert_{\mathfrak p}$ 
an {\it aks}-{\it representation} ({\it associated with} $G/K$).  
Denote by ${\rm ad}_{\mathfrak g}$ the adjoint representation of $\mathfrak g$.  
Let $\mathfrak a_s$ be a maximal split abelian subspace of 
$\mathfrak p$ (see [R] or [OS] about the definition of a maximal split abelian subspace) and 
$\mathfrak p=\mathfrak p_0+\sum\limits_{\alpha\in \triangle_+}\mathfrak p_{\alpha}$ the root 
space decomposition with respect to $\mathfrak a_s$ (i.e., the simultaneously eigenspace 
decomposition of ${\rm ad}_{\mathfrak g}(a)^2$'s ($a\in\mathfrak a_s$)), where the space 
$\mathfrak p_{\alpha}$ is defined by 
$\mathfrak p_{\alpha}:=\{X\in\mathfrak p\,\vert\,{\rm ad}_{\mathfrak g}(a)^2(X)
=\alpha(a)^2X\,\,{\rm for}\,\,{\rm all}\,\,a\in \mathfrak a_s\}$ 
($\alpha\in\mathfrak a_s^{\ast}$) and 
$\triangle_+$ is the positive root system of 
the root system 
$\triangle:=\{\alpha\in\mathfrak a_s^{\ast}\,\vert\,
\mathfrak p_{\alpha}\not=\{0\}\}$ under some lexicographic ordering of 
$\mathfrak a_s^{\ast}$.  
Set $\mathfrak a:=
\mathfrak p_0\,(\supset\mathfrak a_s),\,j:=J_{eK}$ and 
$\langle\,\,,\,\,\rangle_0:=\langle\,\,,\,\,\rangle_{eK}$.  It is shown 
that $\langle\,\,,\,\,\rangle_0\vert_{\mathfrak a_s\times\mathfrak a_s}$ is 
positive (or negative) definite, $\mathfrak a=\mathfrak a_s\oplus 
j\mathfrak a_s$ and $\langle\,\,,\,\,\rangle_0\vert_{\mathfrak a_s\times 
j\mathfrak a_s}=0$.  Note that $\mathfrak p_{\alpha}=\{X\in
\mathfrak p\,\vert\,{\rm ad}_{\mathfrak g}(a)^2(X)=\alpha^{\bf c}(a)^2X\,\,{\rm for}\,\,
{\rm all}\,\,a\in\mathfrak a\}$ holds for each $\alpha\in\triangle_+$, where 
$\alpha^{\bf c}$ is the complexification of $\alpha:\mathfrak a_s\to{\bf R}$ 
(which is a complex linear function over $\mathfrak a_s^{\bf c}=\mathfrak a$) and 
$\alpha^{\bf c}(a)^2X$ means ${\rm Re}(\alpha^{\bf c}(a)^2)X
+{\rm Im}(\alpha^{\bf c}(a)^2)jX$.  Let 
${\it l}_{\alpha}:=(\alpha^{\bf c})^{-1}(0)$ ($\alpha\in\triangle$) and 
$D:=\mathfrak a\setminus\displaystyle{
\mathop{\cup}_{\alpha\in\triangle_+}{\it l}_{\alpha}}$.  Elements of $D$ are said to be 
{\it regular}.  Take $x\in D$ and let 
$M$ be the orbit of the aks-representation ${\rm Ad}_G\vert_{\mathfrak p}$ through $x$.  
From $x\in D$, $M$ is a principal orbit of this representation.  
Denote by $A$ the shape tensor of $M$.  
Take $v\in T^{\perp}_xM(=\mathfrak a)$.  Then we have 
$T_xM=\sum\limits_{\alpha\in\triangle_+}\mathfrak p_{\alpha}$ and 
$A_v\vert_{\mathfrak p_{\alpha}}=-\frac{\alpha^{\bf c}(v)}{\alpha^{\bf c}(x)}{\rm id}$ 
($\alpha\in\triangle_+$).  
From this fact, we see that $M$ is an anti-Kaehler Fredholm submanifold with 
$J$-diagonalizable shape operators.  
Let $\widetilde v$ be the parallel normal vector field of $M$ with $\widetilde v_x=v$.  Then we can show that $A_{\widetilde v_{\rho(k)(x)}}
\vert_{\rho(k)_{\ast x}(\mathfrak p_{\alpha})}=-\frac{\alpha^{\bf c}(v)}{\alpha^{\bf c}(x)}
{\rm id}$ for any $k\in K$.  Hence $M$ is an anti-Kaehler isoparametric submanifold 
with $J$-diagonalizable shape operators.  

\section{Regulalizability of an anti-Kaehler Fredholm submanifold} 
In this section, we shall define the regularizability of an anti-Kaehler Fredholm 
submanifold with $J$-diagonalizable shape operators.  
Let $(M,\langle\,\,,\,\,\rangle_M,J)$ be an anti-Kaehler Fredholm submanifold 
with $J$-diagonalizable shape operators in an infinite dimensional anti-Kaehler space 
$(V,\langle\,\,,\,\,\rangle,\widetilde J)$.  Denote by $A$ the shape tensor of $M$.  
Fix $v\in T^{\perp}M$.  Let $\{\mu_i\,\vert\,i=1,2,\cdots\}$ 
("$\vert\mu_i\vert\,>\,\vert\mu_{i+1}\vert$" or 
"$\vert\mu_i\vert=\vert\mu_{i+1}\vert\,\,\&\,\,
{\rm Re}\,\mu_i>{\rm Re}\,\mu_{i+1}$" or 
"$\vert\mu_i\vert=\vert\mu_{i+1}\vert\,\,\&\,\,
{\rm Re}\,\mu_i={\rm Re}\,\mu_{i+1}\,\,\&\,\,
{\rm Im}\,\mu_i=-{\rm Im}\,\mu_{i+1}>0$") be the set of all $J$-eigenvalues of 
$A_v$ other than zero and $m_i$ the multiplicity of $\mu_i$.  
Then we define the regularized trace 
${\rm Tr}_rA_v$ of $A_v$ by ${\rm Tr}_rA_v:=\sum\limits_{i}m_i\mu_i$.  
Also, we define the trace ${\rm Tr}_{\rm abs}A_v^2$ by 
${\rm Tr}_{\rm abs}A_v^2:=\sum\limits_{i}m_i\vert\mu_i\vert^2$.  
If there exist ${\rm Tr}_rA_v$ and ${\rm Tr}_{\rm abs}A_v^2$ 
for each $v\in T^{\perp}M$, then we say that $M$ is {\it regularizable}.  
It is shown that, if $\mu$ is a $J$-eigenvalue of $A_v$ with multiplicity 
$m$, then so is also the conjugate $\bar{\mu}$ of $\mu$.  Hence we have 
${\rm Tr}_rA_v\in{\Bbb R}$.  Define $H_x\in T^{\perp}_xM$ by 
$\langle H_x,v\rangle={\rm Tr}_rA_v$ ($\forall\,v\in T^{\perp}_xM$).  
We call the normal vector field $H\,(:x\mapsto H_x)$ of $M$ the 
{\it regularized mean curvature vector} of $M$.  

\section{Proof of Theorem A} 
In this section, we shall prove Theorem A.  
For its purpose, we shall prepare some lemmas (and theorems).  
First we shall recall the generalized Chow's theorem, which was proved in 
[HL2].  Let $N$ be a (connected) Hilbert manifold and ${\cal D}$ a set of 
local (smooth) vector fields which are defined over open sets of $N$.  
If two points $x$ 
and $y$ of $N$ can be connected by a piecwise smooth curve each of whose 
smooth segments is an integral curve of a local smooth vector field belonging 
to ${\cal D}$, then we say that $x$ and $y$ are ${\cal D}$-{\it equivalent} 
and we denote this fact by 
$\displaystyle{x\mathop{\sim}_{\cal D}y}$.  Let $\Omega_{\cal D}(x):=\{y\in N
\,\vert\,\displaystyle{y\mathop{\sim}_{\cal D}x\}}$.  The set 
$\Omega_{\cal D}(x)$ is called the {\it set of reachable points of} ${\cal D}$ 
{\it starting from} $x$.  Let ${\cal D}^{\ast}$ be the minimal set 
consisting of local smooth vector fields on open sets of $N$ which satisfies 
the following condition:

\vspace{0.2truecm}

${\cal D}\subset{\cal D}^{\ast}$ and ${\cal D}^{\ast}$ contains the zero 
vector field and, for any $X,Y\in{\cal D}^{\ast}$ and any $a,b\in{\bf R}$, 

$aX+bY$ and $[X,Y]$ (which are defined on the intersection of the domains 
of $X$ 

and $Y$) also belong to ${\cal D}^{\ast}$.  

\vspace{0.2truecm}

\noindent
For each $x\in N$, set 
${\cal D}^{\ast}(x):=\{X_x\,\vert\,X\in{\cal D}^{\ast}\,\,{\rm s.t.}\,\,
x\in{\rm Dom}(X)\}$.  Then the following generalized Chow's theorem holds.  

\vspace{0.2truecm}

\noindent
{\bf Theorem 4.1([HL2])} {\sl If $\overline{{\cal D}^{\ast}(x)}=T_xN$ for each 
$x\in N$, $\overline{\Omega_{\cal D}(x)}=N$ holds for each $x\in N$, where 
$\overline{(\cdot)}$ implies the closure of $(\cdot)$.}

\vspace{0.2truecm}

Let $M$ be as in the statement of Theorem A.  
Denote by $(\langle\,\,,\,\,\rangle_M,J)$ and $A$ the anti-Kaehler structure and the 
shape tensor of $M$, respectively.  For simplicity, we denote 
$\langle\,\,,\,\,\rangle_M$ by $\langle\,\,,\,\,\rangle$.  
Let $\{E_0\}\cup\{E_i\,\vert\,i\in I\}$ the set of all $J$-curvature 
distributions of $M$, where $E_0$ is defined by $(E_0)_x
:=\displaystyle{\mathop{\cap}_{v\in T^{\perp}_xM}
{\rm Ker}\,A_v}$ ($x\in M$).  Also, let $\lambda_i$ and $n_i$ be the $J$-principal 
curvature and the $J$-curvature normal corresponding to $E_i$, respectively.  
Denote by ${\it l}_i^x$ the complex focal hyperplane $(\lambda_i)_x^{-1}(1)$ 
of $M$ at $x$.  Also set $({\it l}_i^x)':=(\lambda_i)_x^{-1}(0)$.  
Fix $x_0\in M$.  For simplicity, set ${\it l}_i:={\it l}_i^{x_0}$ and 
${\it l}'_i:=({\it l}_i^{x_0})'$.  
Let $Q(x_0)$ be the set of all points of $M$ connected with $x_0$ by 
a piecewise smooth curve in $M$ each of whose smooth segments is contained in some 
complex curvature sphere (which may depend on the smooth segment).  
By using the above generalized Chow's theorem, we shall show the following result.  

\vspace{0.2truecm}

\noindent
{\bf Proposition 4.2.} {\sl The set $Q(x_0)$ is dense in $M$.}

\vspace{0.2truecm}

\noindent
{\it Proof.} Let ${\cal D}_E$ be the set of all local (smooth) tangent vector 
fields on open sets of $M$ which is tangent to some $E_i$ ($i\not=0$) at each point of 
the domain.  Define $\Omega_{{\cal D}_E}(x_0),\,{\cal D}^{\ast}_E$ and 
${\cal D}^{\ast}_E(x_0)$ as above.  By imitating the proof of Proposition 5.8 of [HL2], 
it is shown that $\overline{{\cal D}^{\ast}_E(x)}=T_xM$ for each 
$x\in M$.  Hence, $\overline{\Omega_{{\cal D}_E}(x_0)}=M$ follows from Theorem 4.1.  
It is clear that $\Omega_{{\cal D}_E}(x_0)=Q(x_0)$.  Therefore we obtain 
$\overline{Q(x_0)}=M$.  \hspace{3.4truecm}q.e.d.

\vspace{0.2truecm}

For each complex affine subspace $P$ of $T^{\perp}_{x_0}M$, define $I_P$ by 
$$I_P:=\left\{
\begin{array}{ll}
\displaystyle{\{i\in I\,\vert\,(n_i)_{x_0}\in P\}} & 
\displaystyle{(0\notin P)}\\
\displaystyle{\{i\in I\,\vert\,(n_i)_{x_0}\in P\}\cup\{0\}} & 
\displaystyle{(0\in P).}
\end{array}\right.$$
Define a distribution $D_P$ on $M$ by 
$D_P:=\displaystyle{\mathop{\oplus}_{i\in I_P}E_i}$.  

\vspace{0.3truecm}

\noindent
{\bf Lemma 4.3.} {\sl The following statements hold:

{\rm (i)} $M$ is regularizable.

{\rm (ii)} If $0\notin P$, then $I_P$ is finite and 
$\displaystyle{(\mathop{\cap}_{i\in I_P}{\it l}_i)\setminus
(\mathop{\cup}_{i\in I\setminus I_P}{\it l}_i)\not=\emptyset}$.  

{\rm (iii)} If $0\in P$, then $I_P$ is infinite or $I_P=\{0\}$ and 
$\displaystyle{(\mathop{\cap}_{i\in I_P\setminus\{0\}}{\it l}'_i)\setminus
(\mathop{\cup}_{i\in I\setminus I_P}{\it l}'_i)\not=\emptyset}$, where 
$\displaystyle{\mathop{\cap}_{i\in I_P\setminus\{0\}}{\it l}'_i}$ means 
$T^{\perp}_{x_0}M$ when $I_P=\{0\}$.}

\vspace{0.2truecm}

\noindent
{\it Proof.} 
From the discreteness of the complex Coxeter group associated with $M$, 
we can show that $B:=\{(n_i)_{x_0}\,\vert\,i\in I\}$ is described as 
$B=\{\frac{1}{1+a_ij}(n_i)_{x_0}\,\vert\,i\in I_0,\,\,j\in{\Bbb Z}\}$ in terms of some 
finite subset $I_0$ of $I$ and some set $\{a_i\,\vert\,i\in I_0\}$ of complex numbers.  
From this fact, the statements in this lemma follow.  
\hspace{7.7truecm}q.e.d.

\vspace{0.3truecm}

Assume that $0\notin P$. Take 
$v\in\displaystyle{(\mathop{\cap}_{i\in I_P}{\it l}_i)
\setminus(\mathop{\cup}_{i\in I\setminus I_P}{\it l}_i)}$.  
Let $\widetilde v$ be a parallel normal vector field on $M$ 
with $\widetilde v_{x_0}=v$.  
This normal vector field $\widetilde v$ is a focal normal vector field of $M$.  
Let $f_{\widetilde v}$ be the focal map (i.e., the end point map) for $\widetilde v$ 
and $F_{\widetilde v}$ the focal submanifold for $\widetilde v$ (i.e., 
$F_{\widetilde v}=f_{\widetilde v}(M)$).  Also, let $L^{D_P}_x$ be the leaf of $D_P$ through 
$x\in M$.  Note that $L^{D_P}_x=f_{\widetilde v}^{-1}(f_{\widetilde v}(x))$.  
Now we shall show the following homogeneous slice theorem for $M$.  

\vspace{0.3truecm}

\noindent
{\bf Theorem 4.4.} {\sl If $0\notin P$, then the leaf 
$L^{D_P}_x(\subset T^{\perp}_{f_{\widetilde v}(x)}F_{\widetilde v})$ is 
a principal orbit of the direct sum representation 
of some aks-representations and a trivial representation.}

\vspace{0.3truecm}

We shall recall the notion of an anti-Kaehler holonomy system introduced in [Koi4] 
to prove this theorem.  
Let $(W,J,\langle\,\,,\,\,\rangle)$ be a (finite dimensional) 
anti-Kaehler space and $R\,(\in W^{\ast}\otimes W^{\ast}\otimes 
W^{\ast}\otimes W)$ a curvature-like tensor.  Also, let $SO_{AK}(W)$ be 
the identity component of the group 
$\{B\in GL(W)\,\vert\,B^{\ast}\langle\,\,,\,\,\rangle=\langle\,\,,\,\,
\rangle,\,\,[B,J]=0\}$ and $G$ a connected complex Lie subgroup of $SO_{AK}(W)$.  
We call the triple 
$((W,J,\langle\,\,,\,\,\rangle),R,G)$ an {\it anti-Kaehler holonomy 
system} if the following two conditions hold: 

\vspace{0.2truecm}

(i) $J\circ R(w_1,w_2)=R(Jw_1,w_2)=R(w_1,w_2)\circ J$ for all $w_1,w_2\in W$,

(ii) $R(w_1,w_2)\in{\rm Lie}\,G$ for all $w_1,w_2\in W$.  

\vspace{0.2truecm}

\noindent
Furthermore, if the following condition (iii) holds, then we say that the 
triple is {\it symmetric}:

\vspace{0.2truecm}

(iii) $R(gw_1,gw_2)gw_3=gR(w_1,w_2)w_3$ for all $w_i\in W$ ($i=1,2,3$) and 
all $g\in G$.  

\vspace{0.2truecm}

\noindent
Also, if $G$ is weakly irreducible, then we say that the triple is 
{\it weakly irreducible}, where the weakly irreduciblity of $G$ implies that 
there exists no $G$-invariant non-degenerate subspace $W'$ of $W$ with 
$W'\not=\{0\}$ and $W'\not=W$ (where the non-degeneracy of $W'$ impies that 
$\langle\,\,,\,\,\rangle\vert_{W'\times W'}$ is non-degenerate).  
We [Koi4] proved the following fact for a weakly irreducible symmetric anti-Kaehler 
holonomy system.  

\vspace{0.2truecm}

\noindent
{\bf Lemma 4.4.1.} {\sl For a weakly irreducible symmetric anti-Kaehler 
holonomy system \newline
$((W,J,\langle\,\,,\,\,\rangle),R,G)$ with $R\not=0$, 
the $G$-action on $W$ is equivalent to an aks-representation.}

\vspace{0.2truecm}

By using this lemma, we prove Theorem 4.4.  

\vspace{0.2truecm}

\noindent
{\it Proof of Theorem 4.4.} 
Set $x':=f_{\widetilde v}(x)$.  Denote by $\Psi(x')$ the normal holonomy group of 
$F_{\widetilde v}$ at $x'$ and $\Psi^0(x')$ the identity component of $\Psi(x')$.  
Since ${\rm dim}\,T_{x'}^{\perp}F_{\widetilde v}\,<\,\infty$, 
$\Psi^0(x')$ is a Lie subgroup of $SO_{AK}(T_{x'}^{\perp}F_{\widetilde v})$.  
It is clear that $\Psi^0(x')$ is not trivial.  
For simplicity, set $W:=T^{\perp}_{x'}F_{\widetilde v}$.  
Let $W=W_0\oplus W_1\oplus\cdots\oplus W_k$ be the weakly irreducible decomposition of 
the $\Psi^0(x')$-module $W$, where $\Psi^0(x')\vert_{W_0}=\{{\rm id}_{W_0}\}$ and 
$W_i$ ($i=1,\cdots,k$) are (non-trivial) weakly irreducible $\Psi^0(x')$-submodules of $W$.  
For simplicity, set $\Psi^0_i(x'):=\Psi^0(x')\vert_{W_i}$ 
($i=1,\cdots,k$).  Denote by $\widehat A$ the shape tensor of $F_{\widetilde v}$ and 
$R^{\perp}$ the curvature tensor of the normal connection of $F_{\widetilde v}$.  
Also, denote by ${\cal L}^2$ the space of Hilbert-Schmidt operators of 
the Hilbert space $(T_{x'}F_{\widetilde v},\langle\,\,,\,\,\rangle_{V_{\pm}}
\vert_{T_{x'}F_{\widetilde v}\times T_{x'}F_{\widetilde v}})$ and 
$\langle\,\,,\,\,\rangle_{{\cal L}^2}$ the Hilbert-Schmidt inner product of ${\cal L}^2$.  
Define ${\cal R}_i^{\perp}\in W_i^{\ast}\otimes W_i^{\ast}\otimes W_i^{\ast}\otimes W_i$ by 
$$\langle{\cal R}_i^{\perp}(w_1,w_2)w_3,w_4\rangle
:=-\frac12\langle[\widehat A_{w_1},\widehat A_{w_2}],
[\widehat A_{w_3},\widehat A_{w_4}]\rangle_{{\cal L}^2}\quad(w_1,\cdots w_4\in W_i).$$
Here we note that $\widehat A_{w_j}$'s ($j=1,\cdots,4$) are Hilbert-Schmidt operators 
because $M$ (hence $F_{\widetilde v}$) is a regulalizable anti-Kaehler Fredholm submanifold 
with $J$-diagonalizable shape operators.  
From the Ricci equation, $[\widehat A_{w_j},J]=0$ and 
$R^{\perp}(JX,JY)=-R^{\perp}(X,Y)\,\,\,(X,Y\in T_{x'}F_{\widetilde v})$, we can show 
$$\langle{\cal R}_i^{\perp}(w_1,w_2)w_3,w_4\rangle
=2\sum_{j\in{\Bbb N}}\langle R^{\perp}(\widehat A_{w_1}e_j,\widehat A_{w_2}e_j)w_3,w_4
\rangle_{V_{\pm}}\quad(w_1,\cdots w_4\in W_i),$$
where $\{e_j\}_{j=1}^{\infty}$ is a $J$-orthonormal base of $T_{x'}F_{\widetilde v}$.  
By using this relation, we can show that $(W_i,{\cal R}_i^{\perp},\Psi^0_i(x'))$ is 
a weakly irreducible symmetric anti-Kaehler holonomy system.  Also, from 
$R^{\perp}\vert_{T_{x'}F_{\widetilde v}\times T_{x'}F_{\widetilde v}\times W_i}\not=0$, 
we can show ${\cal R}_i^{\perp}\not=0$.  
Hence it follows from Lemma 4.4.1 that 
the $\Psi^0_i(x')$-action on $W_i$ is equivalent to an aks-representation.  
Also, the $\Psi^0_0(x')$-action on $W_0$ is trivial.  
Therefore, since $L^{D_P}_x$ is a principal orbit of $\Psi^0(x')$-action, 
the statement of Theorem 4.4 follows.  \hspace{8.15truecm}q.e.d.

\vspace{0.2truecm}

Set $(W_P)_x:=x+(D_P)_x\oplus{\rm Span}_{\bf C}\{(n_i)_x\,\vert\,i\in I_P
\setminus\{0\}\}$ ($x\in M$).  Let $\gamma:[0,1]\to M$ be a piecewise smooth curve.  
In the sequel, we assume that the domains of all piecewise smooth curves are equal to 
$[0,1]$.  If $\dot{\gamma}(t)\perp(D_P)_{\gamma(t)}$ for each $t\in[0,1]$, then $\gamma$ 
is said to be {\it horizontal with respect to} $D_P$ (or $D_P$-{\it horizontal}).  Let 
$\beta_i$ ($i=1,2$) be curves in $M$.  If 
$L^{D_P}_{\beta_1(t)}=L^{D_P}_{\beta_2(t)}$ for each $t\in[0,1]$, then 
$\beta_1$ and $\beta_2$ are said to be {\it parallel with respect to} $D_P$.  
By imitating the proof of Proposition 1.1 in [HL2], we can show the following fact.  

\vspace{0.2truecm}

\noindent
{\bf Lemma 4.5.} {\sl For each $D_P$-horizontal curve $\gamma$, there exists 
an one-parameter family $\{h^{D_P}_{\gamma,t}\,\vert\,0\leq t\leq 1\}$ of 
holomorphic isometries $h^{D_P}_{\gamma,t}:(W_P)_{\gamma(0)}\to
(W_P)_{\gamma(t)}$ satisfying the following conditions:

{\rm(i)} $h^{D_P}_{\gamma,t}(L^{D_P}_{\gamma(0)})=L^{D_P}_{\gamma(t)}$ 
($0\leq t\leq 1$), 

{\rm(ii)} for any $x\in L^{D_P}_{\gamma(0)}$, 
$t\mapsto h_{\gamma,t}^{D_P}(x)$ is a $D_P$-horizontal curve parallel to 
$\gamma$, 

{\rm(iii)} for any $x\in L^{D_P}_{\gamma(0)}$ and any $i\in I_P$, 
$(h^{D_P}_{\gamma,t})_{\ast x}((E_i)_x)=(E_i)_{h^{D_P}_{\gamma,t}(x)}$.}

\vspace{0.2truecm}

\noindent
{\it Proof.} First we consider the case of $0\notin P$.  Take 
$v\in\displaystyle{\mathop{\cap}_{i\in I_P}{\it l}_i\setminus
(\mathop{\cup}_{i\in I\setminus I_P}{\it l}_i)}$.  Let $\widetilde v$ be the 
parallel normal vector field of $M$ with 
$\widetilde v_{x_0}=v$.  Let 
$\overline{\gamma}:=f_{\widetilde v}\circ\gamma$.  Define a map 
$h_t:(W_P)_{\gamma(0)}\to V$ by 
$h_t(x):=\overline{\gamma}(t)+\bar{\tau}^{\perp}_{\overline{\gamma}\vert_{[0,t]}}
(\overrightarrow{\bar{\gamma}(0)x})$ ($x\in(W_P)_{\gamma(0)}$) (see Figure 1), 
where $\bar{\tau}^{\perp}_{\overline{\gamma}}$ is the parallel translation along 
$\overline{\gamma}$ with respect to the normal connection of 
$F_{\widetilde v}$.  Then it is shown that $\{h_t\,\vert\,0\leq t\leq1\}$ is 
the desired one-parameter family.  Next we consider the case of $0\in P$.  Take 
$v\in\displaystyle{\mathop{\cap}_{i\in I_P\setminus\{0\}}{\it l}'_i
\setminus(\mathop{\cup}_{i\in I\setminus I_P}{\it l}'_i)}$.  
Let $\widetilde v$ be the parallel normal vector field of 
$M$ with $\widetilde v_{x_0}=v$.  We define a map $\nu:M\to S^{\infty}(1)$ 
by $\nu(x):=\widetilde v_x$ ($x\in M$), where $S^{\infty}(1)$ is the unit hypersphere of 
$V$ centered $0$.  Then we have $\nu_{\ast x}=-A_{\widetilde v_x}$ 
($x\in M$).  If $i\in I_P$, then we have 
$\nu_{\ast x}((E_i)_x)=\{-\langle(n_i)_x,\widetilde v_x\rangle X\,\vert\,
X\in(E_i)_x\}=\{0\}$ 
and, if $i\notin I_P$, then we have 
$\nu_{\ast x}((E_i)_x)=\{-\langle(n_i)_x,\widetilde v_x\rangle X\,\vert\,
X\in(E_i)_x\}=(E_i)_x$.  
Hence we have ${\rm Ker}\,\nu_{\ast x}=(D_P)_x$.  Therefore $D_P$ is 
integrable and it gives a foliation on $M$.  Denote by 
$\mathfrak F_P$ this foliation and $D_P^{\perp}$ the orthogonal complementary 
distribution of $\mathfrak F_P$.  Let $U$ be a neighborhood of $\gamma(0)$ in 
$L^{D_P}_{\gamma(0)}$ such that there exists 
a family $\{\psi_t:U\to U_t\,\vert\,0\leq t\leq1\}$ of diffeomophisms 
such that, for any $x\in U$, the curve 
$\displaystyle{\gamma_x\,(\mathop{\Leftrightarrow}_{\rm def}\,
\gamma_x(t):=\psi_t(x))}$ is a $D_P$-horizontal curve, where $U_t$ is a 
neighborhood of $\gamma(t)$ in $L^{D_P}_{\gamma(t)}$.  Note that such a family 
of diffeomorphisms is called an element of holonomy along $\gamma$ 
(with respect to ${\cal F}_P$ and $D_P^{\perp}$) in [BH].  
Let $\triangle$ be a fundamental domain containing $x_0$ of the complex 
Coxeter group of $M$ at $x_0$.  Denote by 
$\triangle_x$ a domain of $T^{\perp}_xM$ given by parallel 
translating $\triangle$ with respect to the normal connection of $M$.  Set 
$\widetilde U:=\displaystyle{\mathop{\cup}_{x\in U}
({\rm Span}_{\bf C}\{(n_i)_x\,\vert\,i\in I_P\setminus\{0\}\}\cap\triangle_x)}$, 
which is an open subset of the affine subspace $(W_P)_{\gamma(0)}$.  
Define a map $h_t:\widetilde U\to(W_P)_{\gamma(t)}$ ($0\leq t\leq1$) by 
$h_t(x+w)=\gamma_x(t)+\tau^{\perp}_{\gamma_x\vert_{[0,t]}}(w)$ 
($x\in U,\,w\in {\rm Span}\{(n_i)_x\,\vert\,i\in I_P\setminus\{0\}\}
\cap\triangle_x$) (see Figure 2).  
By imitating the proof of Lemma 1.2 in [HL2], it is shown that $h_t$ is a 
holomorphic isometry into $(W_P)_{\gamma(t)}$ .  Hence $h_t$ extends to a 
holomorphic isometry of $(W_P)_{\gamma(0)}$ onto $(W_P)_{\gamma(t)}$.  
Denote by $\widetilde h_t$ this holomorphic extension.  
It is shown that $\widetilde h_t$'s gives the desired one-parameter family 
by imitating the discussion in Step 3 of the proof of Proposition 1.1 
in [HL2].  \hspace{12.45truecm}q.e.d.

\newpage


\centerline{
\unitlength 0.1in
\begin{picture}(108.7400, 31.1900)(-49.7000,-38.9900)
%
\special{pn 8}%
\special{ar 4410 3000 2650 970  3.7000526 5.0770201}%
%
\special{pn 8}%
\special{ar 4410 3200 3890 1940  4.0008946 4.9897782}%
%
\special{pn 8}%
\special{ar 4400 3000 3360 1340  3.9041920 5.0201565}%
%
\special{pn 8}%
\special{pa 5320 1700}%
\special{pa 5318 1732}%
\special{pa 5316 1764}%
\special{pa 5310 1796}%
\special{pa 5304 1828}%
\special{pa 5296 1858}%
\special{pa 5286 1888}%
\special{pa 5276 1920}%
\special{pa 5264 1948}%
\special{pa 5250 1978}%
\special{pa 5232 2004}%
\special{pa 5214 2030}%
\special{pa 5190 2052}%
\special{pa 5162 2070}%
\special{pa 5132 2080}%
\special{pa 5100 2076}%
\special{pa 5072 2064}%
\special{pa 5048 2042}%
\special{pa 5028 2018}%
\special{pa 5012 1990}%
\special{pa 4998 1960}%
\special{pa 4988 1930}%
\special{pa 4978 1900}%
\special{pa 4972 1868}%
\special{pa 4968 1836}%
\special{pa 4964 1804}%
\special{pa 4960 1772}%
\special{pa 4960 1740}%
\special{pa 4960 1708}%
\special{pa 4962 1676}%
\special{pa 4964 1644}%
\special{pa 4968 1614}%
\special{pa 4974 1582}%
\special{pa 4980 1550}%
\special{pa 4988 1520}%
\special{pa 4996 1488}%
\special{pa 5008 1458}%
\special{pa 5020 1428}%
\special{pa 5034 1400}%
\special{pa 5050 1372}%
\special{pa 5070 1348}%
\special{pa 5094 1326}%
\special{pa 5122 1310}%
\special{pa 5152 1300}%
\special{pa 5184 1304}%
\special{pa 5212 1318}%
\special{pa 5236 1340}%
\special{pa 5256 1364}%
\special{pa 5272 1392}%
\special{pa 5284 1422}%
\special{pa 5296 1452}%
\special{pa 5304 1484}%
\special{pa 5310 1514}%
\special{pa 5316 1546}%
\special{pa 5320 1578}%
\special{pa 5322 1610}%
\special{pa 5322 1642}%
\special{pa 5322 1674}%
\special{pa 5320 1700}%
\special{sp}%
%
\special{pn 8}%
\special{pa 2630 1900}%
\special{pa 2638 1932}%
\special{pa 2642 1964}%
\special{pa 2648 1994}%
\special{pa 2650 2026}%
\special{pa 2652 2058}%
\special{pa 2652 2090}%
\special{pa 2650 2122}%
\special{pa 2648 2154}%
\special{pa 2642 2186}%
\special{pa 2634 2218}%
\special{pa 2622 2246}%
\special{pa 2606 2274}%
\special{pa 2584 2298}%
\special{pa 2558 2316}%
\special{pa 2526 2324}%
\special{pa 2496 2318}%
\special{pa 2466 2306}%
\special{pa 2440 2286}%
\special{pa 2416 2266}%
\special{pa 2396 2240}%
\special{pa 2376 2216}%
\special{pa 2360 2188}%
\special{pa 2344 2160}%
\special{pa 2330 2130}%
\special{pa 2318 2102}%
\special{pa 2306 2072}%
\special{pa 2296 2042}%
\special{pa 2288 2010}%
\special{pa 2280 1980}%
\special{pa 2272 1948}%
\special{pa 2268 1916}%
\special{pa 2264 1884}%
\special{pa 2260 1852}%
\special{pa 2260 1820}%
\special{pa 2260 1788}%
\special{pa 2260 1756}%
\special{pa 2264 1726}%
\special{pa 2270 1694}%
\special{pa 2278 1662}%
\special{pa 2290 1634}%
\special{pa 2306 1606}%
\special{pa 2328 1582}%
\special{pa 2356 1566}%
\special{pa 2388 1560}%
\special{pa 2418 1566}%
\special{pa 2448 1580}%
\special{pa 2474 1598}%
\special{pa 2496 1620}%
\special{pa 2516 1646}%
\special{pa 2536 1670}%
\special{pa 2552 1698}%
\special{pa 2568 1726}%
\special{pa 2582 1756}%
\special{pa 2594 1784}%
\special{pa 2606 1814}%
\special{pa 2616 1844}%
\special{pa 2624 1876}%
\special{pa 2630 1900}%
\special{sp}%
%
\special{pn 8}%
\special{pa 4000 1670}%
\special{pa 4002 1702}%
\special{pa 4002 1734}%
\special{pa 4000 1766}%
\special{pa 3996 1798}%
\special{pa 3992 1830}%
\special{pa 3988 1862}%
\special{pa 3980 1892}%
\special{pa 3970 1924}%
\special{pa 3960 1954}%
\special{pa 3946 1982}%
\special{pa 3930 2010}%
\special{pa 3910 2034}%
\special{pa 3884 2054}%
\special{pa 3854 2066}%
\special{pa 3824 2068}%
\special{pa 3792 2058}%
\special{pa 3766 2040}%
\special{pa 3744 2016}%
\special{pa 3724 1992}%
\special{pa 3708 1964}%
\special{pa 3694 1934}%
\special{pa 3682 1906}%
\special{pa 3672 1874}%
\special{pa 3664 1844}%
\special{pa 3656 1814}%
\special{pa 3650 1782}%
\special{pa 3646 1750}%
\special{pa 3642 1718}%
\special{pa 3640 1686}%
\special{pa 3640 1654}%
\special{pa 3640 1622}%
\special{pa 3642 1590}%
\special{pa 3644 1558}%
\special{pa 3650 1526}%
\special{pa 3654 1496}%
\special{pa 3662 1464}%
\special{pa 3670 1434}%
\special{pa 3682 1404}%
\special{pa 3696 1374}%
\special{pa 3712 1348}%
\special{pa 3734 1324}%
\special{pa 3758 1304}%
\special{pa 3788 1292}%
\special{pa 3820 1292}%
\special{pa 3850 1302}%
\special{pa 3876 1322}%
\special{pa 3898 1344}%
\special{pa 3918 1370}%
\special{pa 3934 1398}%
\special{pa 3948 1428}%
\special{pa 3960 1456}%
\special{pa 3970 1488}%
\special{pa 3978 1518}%
\special{pa 3984 1550}%
\special{pa 3990 1580}%
\special{pa 3994 1612}%
\special{pa 3998 1644}%
\special{pa 4000 1670}%
\special{sp}%
%
\special{pn 20}%
\special{sh 1}%
\special{ar 2460 1920 10 10 0  6.28318530717959E+0000}%
\special{sh 1}%
\special{ar 2460 1920 10 10 0  6.28318530717959E+0000}%
%
\special{pn 20}%
\special{sh 1}%
\special{ar 2550 1690 10 10 0  6.28318530717959E+0000}%
\special{sh 1}%
\special{ar 2550 1690 10 10 0  6.28318530717959E+0000}%
%
\special{pn 20}%
\special{sh 1}%
\special{ar 3970 1460 10 10 0  6.28318530717959E+0000}%
\special{sh 1}%
\special{ar 3970 1460 10 10 0  6.28318530717959E+0000}%
%
\special{pn 20}%
\special{sh 1}%
\special{ar 2650 2060 10 10 0  6.28318530717959E+0000}%
\special{sh 1}%
\special{ar 2650 2060 10 10 0  6.28318530717959E+0000}%
%
\special{pn 20}%
\special{sh 1}%
\special{ar 3990 1870 10 10 0  6.28318530717959E+0000}%
\special{sh 1}%
\special{ar 3990 1870 10 10 0  6.28318530717959E+0000}%
%
\special{pn 8}%
\special{pa 2450 1920}%
\special{pa 2550 1700}%
\special{fp}%
\special{sh 1}%
\special{pa 2550 1700}%
\special{pa 2504 1752}%
\special{pa 2528 1750}%
\special{pa 2542 1770}%
\special{pa 2550 1700}%
\special{fp}%
%
\special{pn 8}%
\special{pa 3580 950}%
\special{pa 3590 2310}%
\special{pa 4110 2480}%
\special{pa 4110 1080}%
\special{pa 4110 1080}%
\special{pa 3580 950}%
\special{fp}%
%
\special{pn 8}%
\special{pa 2040 1350}%
\special{pa 2330 2620}%
\special{pa 2970 2590}%
\special{pa 2660 1310}%
\special{pa 2660 1310}%
\special{pa 2040 1350}%
\special{fp}%
%
\special{pn 13}%
\special{pa 2460 1910}%
\special{pa 2550 1700}%
\special{fp}%
\special{sh 1}%
\special{pa 2550 1700}%
\special{pa 2506 1754}%
\special{pa 2530 1750}%
\special{pa 2542 1770}%
\special{pa 2550 1700}%
\special{fp}%
%
\special{pn 8}%
\special{pa 3170 3200}%
\special{pa 3800 3300}%
\special{dt 0.045}%
%
\special{pn 8}%
\special{pa 4190 3070}%
\special{pa 3790 3290}%
\special{dt 0.045}%
%
\special{pn 8}%
\special{pa 3800 3300}%
\special{pa 4650 3560}%
\special{dt 0.045}%
\special{sh 1}%
\special{pa 4650 3560}%
\special{pa 4592 3522}%
\special{pa 4600 3544}%
\special{pa 4580 3560}%
\special{pa 4650 3560}%
\special{fp}%
%
\special{pn 8}%
\special{ar 5378 3486 400 160  3.1415927 3.1844498}%
\special{ar 5378 3486 400 160  3.3130212 3.3558784}%
\special{ar 5378 3486 400 160  3.4844498 3.5273069}%
\special{ar 5378 3486 400 160  3.6558784 3.6987355}%
\special{ar 5378 3486 400 160  3.8273069 3.8701641}%
\special{ar 5378 3486 400 160  3.9987355 4.0415927}%
\special{ar 5378 3486 400 160  4.1701641 4.2130212}%
\special{ar 5378 3486 400 160  4.3415927 4.3844498}%
\special{ar 5378 3486 400 160  4.5130212 4.5558784}%
\special{ar 5378 3486 400 160  4.6844498 4.7273069}%
\special{ar 5378 3486 400 160  4.8558784 4.8987355}%
\special{ar 5378 3486 400 160  5.0273069 5.0701641}%
\special{ar 5378 3486 400 160  5.1987355 5.2415927}%
\special{ar 5378 3486 400 160  5.3701641 5.4130212}%
\special{ar 5378 3486 400 160  5.5415927 5.5844498}%
\special{ar 5378 3486 400 160  5.7130212 5.7558784}%
\special{ar 5378 3486 400 160  5.8844498 5.9273069}%
\special{ar 5378 3486 400 160  6.0558784 6.0987355}%
\special{ar 5378 3486 400 160  6.2273069 6.2701641}%
%
\special{pn 8}%
\special{ar 5378 3486 400 150  6.2831853 6.2831853}%
\special{ar 5378 3486 400 150  0.0000000 3.1415927}%
%
\special{pn 8}%
\special{ar 4628 3486 350 540  5.4038175 6.2831853}%
\special{ar 4628 3486 350 540  0.0000000 0.8709035}%
%
\special{pn 8}%
\special{ar 6128 3476 350 540  2.2706892 4.0209605}%
\put(31.1000,-30.2000){\makebox(0,0)[rt]{$L^{D_P}_{\gamma(0)}$}}%
\put(45.8000,-28.3000){\makebox(0,0)[rt]{$L^{D_P}_{\gamma(t)}$}}%
\put(51.3000,-25.6000){\makebox(0,0)[rt]{$\gamma$}}%
%
\special{pn 13}%
\special{ar 4400 3000 3120 1150  4.1173052 4.9986014}%
%
\special{pn 13}%
\special{ar 4400 3000 3380 1340  4.0959929 4.9339661}%
%
\special{pn 20}%
\special{sh 1}%
\special{ar 5290 1910 10 10 0  6.28318530717959E+0000}%
\special{sh 1}%
\special{ar 5290 1910 10 10 0  6.28318530717959E+0000}%
%
\special{pn 20}%
\special{sh 1}%
\special{ar 5140 1700 10 10 0  6.28318530717959E+0000}%
\special{sh 1}%
\special{ar 5140 1700 10 10 0  6.28318530717959E+0000}%
\put(49.5000,-9.1000){\makebox(0,0)[rt]{$\overline{\gamma}$}}%
\put(17.4000,-14.2000){\makebox(0,0)[rt]{$x$}}%
\put(33.5000,-10.9000){\makebox(0,0)[rt]{$h_t(x)$}}%
\put(18.4000,-10.8000){\makebox(0,0)[rt]{$(W_P)_{\gamma(0)}$}}%
\put(33.4000,-7.8000){\makebox(0,0)[rt]{$(W_P)_{\gamma(t)}$}}%
\put(57.0000,-24.6000){\makebox(0,0)[rt]{$F_{\widetilde v}$}}%
\put(15.2000,-17.9000){\makebox(0,0)[rt]{$\overrightarrow{\overline{\gamma}(0)x}$}}%
%
\special{pn 20}%
\special{sh 1}%
\special{ar 3830 1690 10 10 0  6.28318530717959E+0000}%
\special{sh 1}%
\special{ar 3830 1690 10 10 0  6.28318530717959E+0000}%
%
\special{pn 13}%
\special{pa 3820 1690}%
\special{pa 3960 1470}%
\special{fp}%
\special{sh 1}%
\special{pa 3960 1470}%
\special{pa 3908 1516}%
\special{pa 3932 1516}%
\special{pa 3942 1538}%
\special{pa 3960 1470}%
\special{fp}%
\put(42.4000,-34.8000){\makebox(0,0)[rt]{in fact}}%
\put(66.6700,-35.1600){\makebox(0,0)[rt]{{\small non-compact}}}%
\put(53.7000,-11.7000){\makebox(0,0)[lb]{$M$}}%
%
\special{pn 8}%
\special{pa 1760 1510}%
\special{pa 2540 1680}%
\special{dt 0.045}%
\special{sh 1}%
\special{pa 2540 1680}%
\special{pa 2480 1646}%
\special{pa 2488 1670}%
\special{pa 2472 1686}%
\special{pa 2540 1680}%
\special{fp}%
%
\special{pn 8}%
\special{pa 1880 1190}%
\special{pa 2210 1420}%
\special{dt 0.045}%
\special{sh 1}%
\special{pa 2210 1420}%
\special{pa 2168 1366}%
\special{pa 2166 1390}%
\special{pa 2144 1398}%
\special{pa 2210 1420}%
\special{fp}%
%
\special{pn 8}%
\special{pa 1560 1870}%
\special{pa 2490 1820}%
\special{dt 0.045}%
\special{sh 1}%
\special{pa 2490 1820}%
\special{pa 2422 1804}%
\special{pa 2438 1824}%
\special{pa 2426 1844}%
\special{pa 2490 1820}%
\special{fp}%
%
\special{pn 8}%
\special{pa 2860 2920}%
\special{pa 2650 2180}%
\special{dt 0.045}%
\special{sh 1}%
\special{pa 2650 2180}%
\special{pa 2650 2250}%
\special{pa 2666 2232}%
\special{pa 2688 2240}%
\special{pa 2650 2180}%
\special{fp}%
%
\special{pn 8}%
\special{pa 4360 2760}%
\special{pa 3960 1980}%
\special{dt 0.045}%
\special{sh 1}%
\special{pa 3960 1980}%
\special{pa 3974 2048}%
\special{pa 3984 2028}%
\special{pa 4008 2030}%
\special{pa 3960 1980}%
\special{fp}%
%
\special{pn 8}%
\special{pa 5030 2470}%
\special{pa 4710 1870}%
\special{dt 0.045}%
\special{sh 1}%
\special{pa 4710 1870}%
\special{pa 4724 1938}%
\special{pa 4736 1918}%
\special{pa 4760 1920}%
\special{pa 4710 1870}%
\special{fp}%
%
\special{pn 8}%
\special{pa 5580 2420}%
\special{pa 5370 1730}%
\special{dt 0.045}%
\special{sh 1}%
\special{pa 5370 1730}%
\special{pa 5370 1800}%
\special{pa 5386 1782}%
\special{pa 5410 1788}%
\special{pa 5370 1730}%
\special{fp}%
%
\special{pn 8}%
\special{pa 3370 1180}%
\special{pa 3950 1450}%
\special{dt 0.045}%
\special{sh 1}%
\special{pa 3950 1450}%
\special{pa 3898 1404}%
\special{pa 3902 1428}%
\special{pa 3882 1440}%
\special{pa 3950 1450}%
\special{fp}%
%
\special{pn 8}%
\special{pa 3360 860}%
\special{pa 3680 1160}%
\special{dt 0.045}%
\special{sh 1}%
\special{pa 3680 1160}%
\special{pa 3646 1100}%
\special{pa 3642 1124}%
\special{pa 3618 1130}%
\special{pa 3680 1160}%
\special{fp}%
%
\special{pn 8}%
\special{pa 4810 1120}%
\special{pa 4520 1670}%
\special{dt 0.045}%
\special{sh 1}%
\special{pa 4520 1670}%
\special{pa 4570 1620}%
\special{pa 4546 1624}%
\special{pa 4534 1602}%
\special{pa 4520 1670}%
\special{fp}%
%
\special{pn 8}%
\special{pa 5330 1140}%
\special{pa 4880 1420}%
\special{dt 0.045}%
\special{sh 1}%
\special{pa 4880 1420}%
\special{pa 4948 1402}%
\special{pa 4926 1392}%
\special{pa 4926 1368}%
\special{pa 4880 1420}%
\special{fp}%
%
\special{pn 8}%
\special{pa 3170 3200}%
\special{pa 3790 3300}%
\special{fp}%
%
\special{pn 8}%
\special{pa 4200 3070}%
\special{pa 3760 3300}%
\special{fp}%
\special{pa 3760 3300}%
\special{pa 3760 3300}%
\special{fp}%
%
\special{pn 8}%
\special{pa 3770 3290}%
\special{pa 4650 3560}%
\special{fp}%
\special{sh 1}%
\special{pa 4650 3560}%
\special{pa 4592 3522}%
\special{pa 4600 3544}%
\special{pa 4580 3560}%
\special{pa 4650 3560}%
\special{fp}%
%
\special{pn 8}%
\special{pa 2540 1070}%
\special{pa 3864 1612}%
\special{dt 0.045}%
\special{sh 1}%
\special{pa 3864 1612}%
\special{pa 3810 1568}%
\special{pa 3814 1592}%
\special{pa 3794 1604}%
\special{pa 3864 1612}%
\special{fp}%
\put(25.0000,-10.7000){\makebox(0,0)[rb]{$\bar{\tau}_{\bar{\gamma}\vert_{[0,t]}}^{\perp}(\overrightarrow{\bar{\gamma}(0)x})$}}%
%
\special{pn 8}%
\special{pa 1880 2550}%
\special{pa 2440 1950}%
\special{dt 0.045}%
\special{sh 1}%
\special{pa 2440 1950}%
\special{pa 2380 1986}%
\special{pa 2404 1990}%
\special{pa 2410 2012}%
\special{pa 2440 1950}%
\special{fp}%
%
\special{pn 8}%
\special{pa 3070 2360}%
\special{pa 2646 2062}%
\special{dt 0.045}%
\special{sh 1}%
\special{pa 2646 2062}%
\special{pa 2688 2118}%
\special{pa 2690 2094}%
\special{pa 2712 2084}%
\special{pa 2646 2062}%
\special{fp}%
%
\special{pn 8}%
\special{pa 3580 2480}%
\special{pa 3820 1710}%
\special{dt 0.045}%
\special{sh 1}%
\special{pa 3820 1710}%
\special{pa 3782 1768}%
\special{pa 3804 1762}%
\special{pa 3820 1780}%
\special{pa 3820 1710}%
\special{fp}%
%
\special{pn 8}%
\special{pa 4340 2160}%
\special{pa 4000 1890}%
\special{dt 0.045}%
\special{sh 1}%
\special{pa 4000 1890}%
\special{pa 4040 1948}%
\special{pa 4042 1924}%
\special{pa 4066 1916}%
\special{pa 4000 1890}%
\special{fp}%
%
\special{pn 8}%
\special{pa 5600 1470}%
\special{pa 5150 1680}%
\special{dt 0.045}%
\special{sh 1}%
\special{pa 5150 1680}%
\special{pa 5220 1670}%
\special{pa 5198 1658}%
\special{pa 5202 1634}%
\special{pa 5150 1680}%
\special{fp}%
%
\special{pn 8}%
\special{pa 5730 2190}%
\special{pa 5320 1930}%
\special{dt 0.045}%
\special{sh 1}%
\special{pa 5320 1930}%
\special{pa 5366 1984}%
\special{pa 5366 1960}%
\special{pa 5388 1950}%
\special{pa 5320 1930}%
\special{fp}%
\put(20.0000,-26.0000){\makebox(0,0)[rt]{$\bar{\gamma}(0)$}}%
\put(30.0000,-24.0000){\makebox(0,0)[lt]{$\gamma(0)$}}%
\put(37.4000,-24.6000){\makebox(0,0)[rt]{$\bar{\gamma}(t)$}}%
\put(56.3000,-13.6000){\makebox(0,0)[lt]{$\bar{\gamma}(1)$}}%
\put(57.9000,-21.6000){\makebox(0,0)[lt]{$\gamma(1)$}}%
\put(43.2000,-22.2000){\makebox(0,0)[lt]{$\gamma(t)$}}%
\end{picture}%
\hspace{17.8truecm}
}

\vspace{0.9truecm}

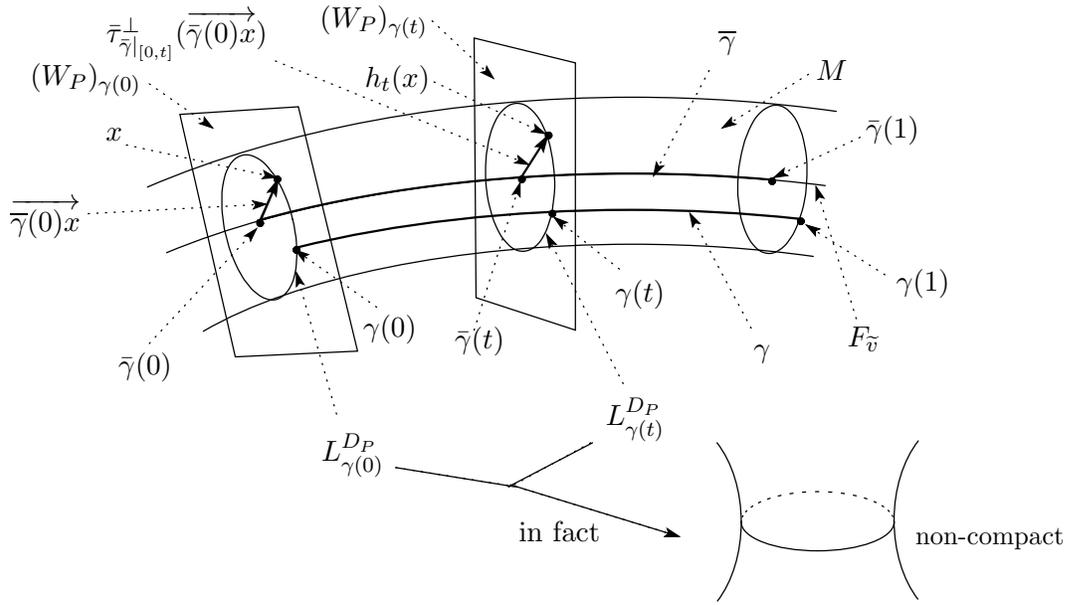
\centerline{{\bf Figure 1.}}


\vspace{1.6truecm}

\centerline{
\unitlength 0.1in
\begin{picture}( 47.6600, 23.1000)(  1.3700,-27.0000)
%
\special{pn 8}%
\special{ar 1314 2000 186 400  0.0000000 6.2831853}%
%
\special{pn 8}%
\special{ar 1314 2000 186 400  0.0000000 6.2831853}%
%
\special{pn 8}%
\special{ar 4104 2000 186 400  0.0000000 6.2831853}%
%
\special{pn 8}%
\special{ar 4104 2000 186 400  0.0000000 6.2831853}%
%
\special{pn 8}%
\special{ar 1872 2000 186 400  0.0000000 6.2831853}%
%
\special{pn 8}%
\special{ar 1872 2000 186 400  0.0000000 6.2831853}%
%
\special{pn 8}%
\special{ar 2988 2000 186 400  0.0000000 6.2831853}%
%
\special{pn 8}%
\special{ar 2988 2000 186 400  0.0000000 6.2831853}%
%
\special{pn 8}%
\special{pa 1314 1600}%
\special{pa 4104 1600}%
\special{fp}%
%
\special{pn 8}%
\special{pa 1314 2400}%
\special{pa 4104 2400}%
\special{fp}%
%
\special{pn 8}%
\special{pa 1314 2000}%
\special{pa 4104 2000}%
\special{dt 0.045}%
%
\special{pn 8}%
\special{pa 1472 1770}%
\special{pa 4262 1770}%
\special{dt 0.045}%
%
\special{pn 8}%
\special{pa 1314 2000}%
\special{pa 1314 800}%
\special{fp}%
\special{pa 1314 800}%
\special{pa 4104 800}%
\special{fp}%
\special{pa 4104 800}%
\special{pa 4104 2000}%
\special{fp}%
%
\special{pn 8}%
\special{pa 1314 2000}%
\special{pa 1908 1000}%
\special{fp}%
%
\special{pn 8}%
\special{pa 4104 2000}%
\special{pa 4698 1000}%
\special{fp}%
%
\special{pn 8}%
\special{pa 1908 1000}%
\special{pa 4688 1000}%
\special{fp}%
%
\special{pn 8}%
\special{pa 1872 2000}%
\special{pa 1872 1600}%
\special{dt 0.045}%
%
\special{pn 8}%
\special{pa 1872 2000}%
\special{pa 2020 1770}%
\special{dt 0.045}%
%
\special{pn 8}%
\special{pa 2988 2000}%
\special{pa 2988 1600}%
\special{dt 0.045}%
%
\special{pn 8}%
\special{pa 2988 2000}%
\special{pa 3136 1770}%
\special{dt 0.045}%
%
\special{pn 13}%
\special{ar 1872 2000 186 400  4.7378792 5.6898352}%
%
\special{pn 13}%
\special{ar 2988 2000 186 400  4.6812464 5.6645800}%
%
\special{pn 13}%
\special{ar 1872 2000 186 400  4.7123890 6.2831853}%
%
\special{pn 13}%
\special{ar 2988 2000 186 400  4.6803487 6.2831853}%
%
\special{pn 13}%
\special{pa 2988 1600}%
\special{pa 2988 1200}%
\special{fp}%
\special{sh 1}%
\special{pa 2988 1200}%
\special{pa 2968 1268}%
\special{pa 2988 1254}%
\special{pa 3008 1268}%
\special{pa 2988 1200}%
\special{fp}%
%
\special{pn 13}%
\special{pa 3146 1770}%
\special{pa 3378 1450}%
\special{fp}%
\special{sh 1}%
\special{pa 3378 1450}%
\special{pa 3324 1492}%
\special{pa 3348 1494}%
\special{pa 3356 1516}%
\special{pa 3378 1450}%
\special{fp}%
%
\special{pn 20}%
\special{sh 1}%
\special{ar 2988 1600 10 10 0  6.28318530717959E+0000}%
\special{sh 1}%
\special{ar 2988 1600 10 10 0  6.28318530717959E+0000}%
%
\special{pn 20}%
\special{sh 1}%
\special{ar 3146 1770 10 10 0  6.28318530717959E+0000}%
\special{sh 1}%
\special{ar 3146 1770 10 10 0  6.28318530717959E+0000}%
%
\special{pn 20}%
\special{sh 1}%
\special{ar 1872 1600 10 10 0  6.28318530717959E+0000}%
\special{sh 1}%
\special{ar 1872 1600 10 10 0  6.28318530717959E+0000}%
%
\special{pn 20}%
\special{sh 1}%
\special{ar 2978 1200 10 10 0  6.28318530717959E+0000}%
\special{sh 1}%
\special{ar 2978 1200 10 10 0  6.28318530717959E+0000}%
%
\special{pn 20}%
\special{sh 1}%
\special{ar 3378 1450 10 10 0  6.28318530717959E+0000}%
\special{sh 1}%
\special{ar 3378 1450 10 10 0  6.28318530717959E+0000}%
%
\special{pn 8}%
\special{pa 4904 1210}%
\special{pa 3396 1440}%
\special{dt 0.045}%
\special{sh 1}%
\special{pa 3396 1440}%
\special{pa 3466 1450}%
\special{pa 3450 1432}%
\special{pa 3460 1410}%
\special{pa 3396 1440}%
\special{fp}%
\put(27.5400,-13.2000){\makebox(0,0)[rt]{$x$}}%
\put(24.3800,-5.5000){\makebox(0,0)[rt]{$x+w$}}%
\put(31.5500,-13.3000){\makebox(0,0)[rt]{$w$}}%
\put(16.2900,-27.0000){\makebox(0,0)[rt]{$\gamma$}}%
\put(24.6600,-26.9000){\makebox(0,0)[rt]{$\gamma_x$}}%
%
\special{pn 20}%
\special{sh 1}%
\special{ar 2030 1770 10 10 0  6.28318530717959E+0000}%
\special{sh 1}%
\special{ar 2030 1770 10 10 0  6.28318530717959E+0000}%
%
\special{pn 8}%
\special{pa 2326 2060}%
\special{pa 2048 1770}%
\special{fp}%
\special{sh 1}%
\special{pa 2048 1770}%
\special{pa 2080 1832}%
\special{pa 2086 1810}%
\special{pa 2110 1804}%
\special{pa 2048 1770}%
\special{fp}%
\put(24.3800,-20.9000){\makebox(0,0)[rt]{$\gamma(t)$}}%
\put(35.2600,-21.0000){\makebox(0,0)[rt]{$\gamma_x(t)$}}%
\put(11.2700,-13.5000){\makebox(0,0)[rt]{$\gamma(0)$}}%
\put(53.7700,-19.6000){\makebox(0,0)[rt]{$\tau^{\perp}_{\gamma_x\vert_{[0,t]}}(w)$}}%
\put(55.6200,-11.2000){\makebox(0,0)[rt]{$h_t(x+w)$}}%
\put(44.3700,-3.9000){\makebox(0,0)[rt]{$(W_P)_{\gamma(0)}$}}%
\put(51.3500,-3.9000){\makebox(0,0)[rt]{$(W_P)_{\gamma(t)}$}}%
%
\special{pn 13}%
\special{pa 1314 1600}%
\special{pa 4104 1610}%
\special{fp}%
\put(49.6800,-13.6000){\makebox(0,0)[rt]{$L^{D_P}_{\gamma(0)}$}}%
%
\special{pn 8}%
\special{pa 1146 1430}%
\special{pa 1862 1590}%
\special{dt 0.045}%
\special{sh 1}%
\special{pa 1862 1590}%
\special{pa 1802 1556}%
\special{pa 1810 1578}%
\special{pa 1794 1596}%
\special{pa 1862 1590}%
\special{fp}%
%
\special{pn 8}%
\special{pa 2430 730}%
\special{pa 2970 1180}%
\special{dt 0.045}%
\special{sh 1}%
\special{pa 2970 1180}%
\special{pa 2932 1122}%
\special{pa 2928 1146}%
\special{pa 2906 1154}%
\special{pa 2970 1180}%
\special{fp}%
%
\special{pn 8}%
\special{pa 2764 1420}%
\special{pa 2988 1600}%
\special{dt 0.045}%
\special{sh 1}%
\special{pa 2988 1600}%
\special{pa 2948 1544}%
\special{pa 2946 1568}%
\special{pa 2924 1574}%
\special{pa 2988 1600}%
\special{fp}%
%
\special{pn 8}%
\special{pa 3340 2050}%
\special{pa 3156 1780}%
\special{dt 0.045}%
\special{sh 1}%
\special{pa 3156 1780}%
\special{pa 3176 1846}%
\special{pa 3186 1824}%
\special{pa 3210 1824}%
\special{pa 3156 1780}%
\special{fp}%
%
\special{pn 8}%
\special{pa 4624 2060}%
\special{pa 3238 1660}%
\special{dt 0.045}%
\special{sh 1}%
\special{pa 3238 1660}%
\special{pa 3296 1698}%
\special{pa 3290 1676}%
\special{pa 3308 1660}%
\special{pa 3238 1660}%
\special{fp}%
%
\special{pn 8}%
\special{pa 4234 570}%
\special{pa 3982 900}%
\special{dt 0.045}%
\special{sh 1}%
\special{pa 3982 900}%
\special{pa 4038 860}%
\special{pa 4014 858}%
\special{pa 4008 836}%
\special{pa 3982 900}%
\special{fp}%
%
\special{pn 8}%
\special{pa 4838 570}%
\special{pa 4400 1150}%
\special{dt 0.045}%
\special{sh 1}%
\special{pa 4400 1150}%
\special{pa 4456 1110}%
\special{pa 4432 1108}%
\special{pa 4424 1086}%
\special{pa 4400 1150}%
\special{fp}%
%
\special{pn 8}%
\special{pa 4568 1410}%
\special{pa 3880 1600}%
\special{dt 0.045}%
\special{sh 1}%
\special{pa 3880 1600}%
\special{pa 3950 1602}%
\special{pa 3930 1586}%
\special{pa 3938 1564}%
\special{pa 3880 1600}%
\special{fp}%
%
\special{pn 8}%
\special{pa 2346 2630}%
\special{pa 3156 1900}%
\special{dt 0.045}%
\special{sh 1}%
\special{pa 3156 1900}%
\special{pa 3092 1930}%
\special{pa 3116 1936}%
\special{pa 3120 1960}%
\special{pa 3156 1900}%
\special{fp}%
%
\special{pn 8}%
\special{pa 1556 2650}%
\special{pa 2038 1900}%
\special{dt 0.045}%
\special{sh 1}%
\special{pa 2038 1900}%
\special{pa 1986 1946}%
\special{pa 2010 1946}%
\special{pa 2020 1968}%
\special{pa 2038 1900}%
\special{fp}%
%
\special{pn 8}%
\special{pa 3982 2620}%
\special{pa 3768 2260}%
\special{dt 0.045}%
\special{sh 1}%
\special{pa 3768 2260}%
\special{pa 3786 2328}%
\special{pa 3796 2306}%
\special{pa 3820 2308}%
\special{pa 3768 2260}%
\special{fp}%
\put(39.5000,-26.8000){\makebox(0,0)[lt]{$M$}}%
%
\special{pn 8}%
\special{pa 2990 1610}%
\special{pa 2990 800}%
\special{dt 0.045}%
%
\special{pn 8}%
\special{pa 3150 1770}%
\special{pa 3740 1000}%
\special{dt 0.045}%
\end{picture}%
\hspace{2.5truecm}}

\vspace{0.9truecm}

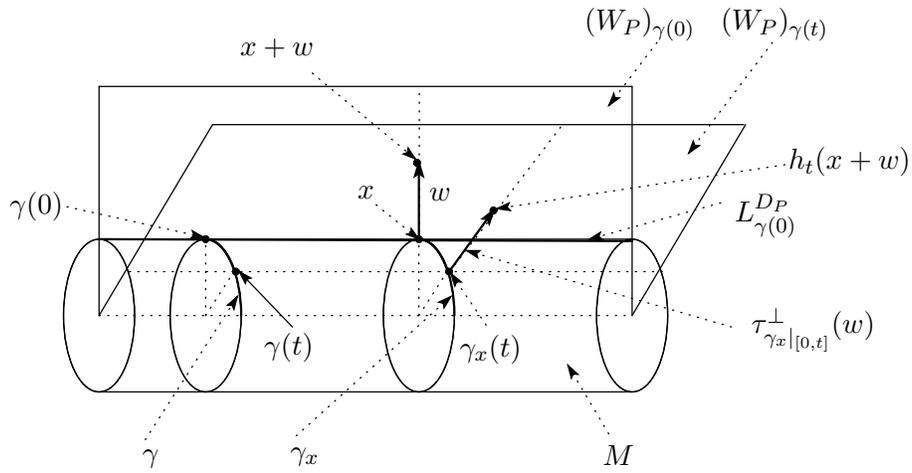
\centerline{{\bf Figure 2.}}

\newpage

\noindent
Fix $x_0\in M$ and $i_0\in I\cup\{0\}$. Take a complex affine subspace $P_{i_0}$ of 
$T^{\perp}_{x_0}M$ with $I_{P_{i_0}}=\{i_0\}$.  Note that 
$D_{P_{i_0}}$ is equal to $E_{i_0}$.  
Denote by $\Phi_{i_0}(x_0)$ the group of holomorphic isometries of 
$(W_{P_{i_0}})_{x_0}$ generated by 
$\{h^{E_{i_0}}_{\gamma,1}\,\vert\,\gamma\,:\,E_{i_0}-{\rm horizontal}
\,\,{\rm curve}\,\,{\rm s.t.}\,\,\gamma(0),\gamma(1)\in L_{x_0}^{E_{i_0}}\}$, 
where $L^{E_{i_0}}_{x_0}$ is the integral manifold of $E_{i_0}$ through 
$x_0$.  Also, denote by $\Phi^0_{i_0}(x_0)$ the identity component of 
$\Phi_{i_0}(x_0)$ and $\Phi^0_{i_0}(x_0)_{x_0}$ the isotropy subgroup of 
$\Phi^0_{i_0}(x_0)$ at $x_0$.  
Define a ${\rm Ad}_{\Phi^0_{i_0}(x_0)}(\Phi^0_{i_0}(x_0))$-invariant 
non-degenerate inner product $\langle\,\,,\,\,\rangle$ of the Lie algebra 
${\rm Lie}\,\Phi^0_{i_0}(x_0)$ of $\Phi^0_{i_0}(x_0)$ by 
$$\langle X,Y\rangle:=B(X,Y)+{\rm Tr}(X\circ Y)\quad(X,Y\in{\rm Lie}\,
\Phi^0_{i_0}(x_0)),$$
where $B$ is the Killing form of ${\rm Lie}\,\Phi^0_{i_0}(x_0)$ and $X\circ Y$ implies 
the composition of $X$ and $Y$ regarded as linear transformations of 
$(W_{P_{i_0}})_{x_0}$.  Take 
$X\in{\rm Lie}\,\Phi^0_{i_0}(x_0)\ominus{\rm Lie}\,\Phi^0_{i_0}(x_0)_{x_0}$.  
Set $g(t):=\exp tX$ and $\gamma(t):=g(t)x_0$, where $\exp$ is the exponential 
map of $\Phi^0_{i_0}(x_0)$.  
It is clear that $\gamma$ is 
an $E_i$-horizontal curve for each $i\in I$ with $i\not=i_0$.  
Let $F_{\gamma}$ be the holomorphic isometry of $V$ satisfying 
$F_{\gamma}(\gamma(0))=\gamma(1)$ and 
$$(F_{\gamma})_{\ast\gamma(0)}=\left\{
\begin{array}{ll}
\displaystyle{g(1)_{\ast\gamma(0)}}&
\displaystyle{{\rm on}\,\,(E_{i_0})_{\gamma(0)}}\\
\displaystyle{(h^{E_i}_{\gamma,1})_{\ast\gamma(0)}}&
\displaystyle{{\rm on}\,\,(E_i)_{\gamma(0)}\,\,(i\in(I\cup\{0\})\setminus\{i_0\})}\\
\displaystyle{\tau^{\perp}_{\gamma}}&
\displaystyle{{\rm on}\,\,T^{\perp}_{\gamma(0)}M.}
\end{array}
\right.$$
In similar to Theorem 4.1 of [HL2], we have the following fact.  

\vspace{0.3truecm}

\noindent
{\bf Proposition 4.6.} {\sl The holomorphic isometry $F_{\gamma}$ preserves 
$M$ invariantly (i.e., $F_{\gamma}(M)=M$).  Furthermore, it preserves $E_i$ 
($i\in I\cup\{0\}$) invariantly (i.e., $F_{\gamma\ast}(E_i)=E_i$).}

\vspace{0.3truecm}

To show this proposition, we prepare some lemmas.  
By imitating the proof (P163$\sim$166) of Proposition 3.1 in [HL2], we can 
show the following fact.  

\vspace{0.5truecm}

\noindent
{\bf Lemma 4.6.1.} {\sl Let $N$ and $\widehat N$ be full irreducible 
anti-Kaehler isoparametric submanifolds with $J$-diagonalizable shape operators 
in an infinite dimensional anti-Kaehler space.  
If ${\rm codim}_{\bf c}N={\rm codim}_{\bf c}\widehat N\geq 2$, $N\cap\widehat N
\not=\emptyset$ and, for some $x_0\in N\cap\widehat N$, 
$T_{x_0}N=T_{x_0}\widehat N$ and if there exists a complex affine line 
${\it l}_0$ of $T^{\perp}_{x_0}N(=T^{\perp}_{x_0}\widehat N)$ such that 
$L^{D_{\it l}}_{x_0}=L^{{\widehat D}_{\it l}}_{x_0}$ for any complex affine 
line ${\it l}$ of $T^{\perp}_{x_0}N$ with ${\it l}\not={\it l}_0$, then 
$N=\widehat N$ holds, where $D_{\it l}$ (resp. $\widehat D_{\it l}$) is 
the integrable distribution on $N$ (resp. $\widehat N$) defined for ${\it l}$ 
in similar to $D_P$.}

\vspace{0.5truecm}

\noindent
{\it Proof.} Let $\{\lambda_i\,\vert\,i\in I\}$ (resp. 
$\{\widehat{\lambda}_i\,\vert\,i\in\widehat I\}$) be the set of all 
$J$-principal curvatures of $N$ (resp. $\widehat N$), ${\bf n}_i$ (resp. 
$\widehat{\bf n}_i$) the $J$-curvature normal corresponding to 
$\lambda_i$ (resp. $\widehat{\lambda}_i$) and $E_i$ (resp. $\widehat E_i$) 
the $J$-curvature distribution corresponding to $\lambda_i$ (resp. 
$\widehat{\lambda}_i$).  Denote by $A$ (resp. $\widehat A$) the shape tensor 
of $N$ (resp. $\widehat N$).  
Let $E_0$ be the $J$-curvature distribution on $N$ with 
$\displaystyle{(E_0)_x:=\mathop{\cap}_{v\in T_x^{\perp}N}{\rm Ker}\,A_v}$ 
($x\in N$) and $\widehat E_0$ the $J$-curvature distribution on 
$\widehat N$ with $\displaystyle{(\widehat E_0)_x
:=\mathop{\cap}_{v\in T_x^{\perp}\widehat N}{\rm Ker}\,\widehat A_v}$ 
($x\in \widehat N$).  For each $x\in N$ (resp. $\hat x\in\widehat N$), 
let $Q_0(x)$ (resp. $\widehat Q_0(\hat x)$) be the set of all points of $N$ 
(resp. $\widehat N$) connected with $x$ (resp. $\hat x$) by a piecewise smooth curve in $N$ 
(resp. $\widehat N$) each of whose smooth segments is contained in some 
complex curvature sphere in $N$ (resp. $\widehat N$) or some integral manifold of $E_0$ 
(resp. $\widehat E_0$).  Take any $x\in Q_0(x_0)$.  
There exists a sequence $\{x_0,x_1,\cdots,x_k(=x)\}$ such that, for each 
$j\in\{1,\cdots,k\}$, $x_j\in\displaystyle{(\mathop{\cup}_{i\in I}
L^{E_i}_{x_{j-1}})\cup L^{E_0}_{x_{j-1}}}$ holds.  Assume that there exists 
$j_0\in\{1,\cdots,k\}$ such that $x_{j_0}\in L^{E_{i_0}}_{x_{j_0-1}}$ 
for some $i_0\in I$ with $(n_{i_0})_{x_0}\in{\it l}_0$.  Since $N$ is 
irreducible, the complex Coxeter group associated with $N$ is not 
decomposable.  Furthermore, since $N$ is full, the group is not decomposed trivially.  
Hence, according to Lemma 3.8 of [Koi3], we can find 
a $J$-curvature normal $n_{i_1}$ of $N$ satisfying 
$(n_{i_1})_{x_0}\notin{\rm Span}_{\bf C}\{(n_{i_0})_{x_0}\}\cup
{\rm Span}_{\bf C}\{(n_{i_0})_{x_0}\}^{\perp}$ (see the final part of the 
first paragraph of Section 2), where we use also 
${\rm codim}_{\bf c}N\geq2$.  Furthermore, since $n_{i_1}$ is a $J$-curvature 
normal, so are also infinitely many complex-constant-multiples of $n_{i_1}$.  
Hence we may assume that $(n_{i_1})_{x_0}$ does not belong to 
${\it l}_0$ by replacing $n_{i_1}$ to a complex-constant-multiple of $n_{i_1}$ 
if necessary.  
Denote by ${\it l}_{i_0i_1}$ the affine line in $T^{\perp}_{x_0}N$ through 
$(n_{i_0})_{x_0}$ and $(n_{i_1})_{x_0}$, and set 
$D_{i_0i_1}:=D_{{\it l}_{i_0i_1}}$ for simplicity.  
According to Theorem 4.4, 
$L^{D_{{i_0i_1}}}_{x_{j_0-1}}$ is a principal orbit of the direct sum representation 
of some aks-representations and a trivial representation and hence it is 
an anti-Kaehler isoparametric submanifold with $J$-diagonalizable shape 
operators in $(W_{{\it l}_{i_0i_1}})_{x_{j_0-1}}$ of complex codimension two.  
Furthermore, since both $(n_{i_0})_{x_0}$ and $(n_{i_1})_{x_0}$ 
are $J$-curvature normals of $L^{D_{{i_0i_1}}}_{x_{j_0-1}}
(\subset(W_{{\it l}_{i_0i_1}})_{x_{j_0-1}})$ and since they are not 
orthogonal, it follows from Lemma 3.8 of [Koi3] that 
$L^{D_{{i_0i_1}}}_{x_{j_0-1}}$ is irreducible.  
Hence, by the anti-Kaehler version of Theorem D of [HOT], $x_{j_0-1}$ can 
be joined to $x_{j_0}$ by a piecewise smooth curve each of whose smooth 
segments is tangent to one of $E_i$'s ($i\in I\,\,{\rm s.t.}\,\,(n_i)_{x_0}
\in{\it l}_{i_0i_1}\,\,{\rm and}\,\,(n_i)_{x_0}\not=(n_{i_0})_{x_0}$).  
Therefore, we can find a sequence $\{x_0,x'_1,\cdots,x'_{k'}(=x)\}$ 
such that, for each $j\in\{1,\cdots,k'\}$, $x'_j\in\displaystyle{
\left(\mathop{\cup}_{i\in I\,\,{\rm s.t.}\,\,(n_i)_{x_0}\notin{\it l}_0}
L^{E_i}_{x'_{j-1}}\right)\cup L^{E_0}_{x'_{j-1}}}$ holds.  Hence it follows from 
Lemma 4.6.2 (see below) that 
$x'_1\in\widehat Q_0(x_0),\,x'_2\in\widehat Q_0(x'_1),\,\cdots,
x'_{k'-1}\in\widehat Q_0(x'_{k'-2})$ and 
$x\in\widehat Q_0(x'_{k'-1})$ inductively.  Therefore we have 
$x\in\widehat Q_0(x_0)$.  From the arbitrariness of $x$, it follows that 
$Q_0(x_0)\subset\widehat Q_0(x_0)$.  Similarly we can show 
$\widehat Q_0(x_0)\subset Q_0(x_0)$.  Thus we obtain 
$Q_0(x_0)=\widehat Q_0(x_0)$ and hence 
$\overline{Q_0(x_0)}=\overline{\widehat Q_0(x_0)}$.  
Let ${\cal D}_E^0$ (resp. $\widehat{{\cal D}}_E^0$) be the set of all local 
(smooth) vector fields of $N$ (resp. $\widehat N$) which is tangent to some 
$E_i$ (resp. $\widehat E_i$) ($i\in I\cup\{0\}$) at each point of 
the domain.  
Since $\overline{({\cal D}_E^0)^{\ast}(x)}=
\overline{(E_0)_x\oplus(\displaystyle{\mathop{\oplus}_{i\in I}}(E_i)_x)}
=T_xN$ for each $x\in N$, it follows from Theorem 4.1 that 
$\overline{\Omega_{{\cal D}_E^0}(x_0)}=N$.  
Similarly, we have 
$\overline{\Omega_{\widehat{\cal D}_E^0}(x_0)}=\widehat N$.  Also, 
it is clear that $\Omega_{{\cal D}_E^0}(x_0)=Q_0(x_0)$ and 
$\Omega_{\widehat{\cal D}_E^0}(x_0)=\widehat Q_0(x_0)$.  
Therefore we obtain $N=\widehat N$.\hspace{2truecm}q.e.d.

\vspace{0.5truecm}

\noindent
{\bf Lemma 4.6.2.} {\sl Let $N,\,\widehat N,\,x_0$ and ${\it l}_0$ be as in 
Lemma 4.6.1.  Then we have $L^{D_{\it l}}_x=L^{\widehat D_{\it l}}_x$ 
for any $x\in L^{E_0}_{x_0}\cup
(\displaystyle{\mathop{\cup}_{i\in I\,\,{\rm s.t.}\,\,(n_i)_{x_0}\notin
{\it l}_0}L^{E_i}_{x_0})}$ and any complex affine line 
${\it l}$ of $T^{\perp}_{x_0}N$ with ${\it l}\not={\it l}_0$.  Also, 
we have $T_xN=T_x\widehat N$ for any $x\in L^{E_0}_{x_0}\cup
(\displaystyle{\mathop{\cup}_{i\in I\,\,{\rm s.t.}\,\,(n_i)_{x_0}\notin
{\it l}_0}L^{E_i}_{x_0})}$.}

\vspace{0.5truecm}

\noindent
{\it Proof.} Assume that $x\in L^{E_{i_0}}_{x_0}$, where $i_0$ is an element of 
$\{i\in I\,\vert\,(n_i)_{x_0}\notin{\it l}_0\}\cup\{0\}$.  
Take any complex affine line ${\it l}$ of $T^{\perp}_{x_0}N$ with ${\it l}\not={\it l}_0$.  
In case of $(n_{i_0})_{x_0}\in{\it l}$, we have 
$x\in L^{E_{i_0}}_{x_0}\subset L^{D_{\it l}}_{x_0}=L^{{\widehat D}_{\it l}}_{x_0}$ and 
hence $L^{D_{\it l}}_x=L^{{\widehat D}_{\it l}}_x$.  
We consider the case of $(n_{i_0})_{x_0}\notin{\it l}$.  Take a curve 
$\gamma:[0,1]\to L^{E_{i_0}}_{x_0}$ with $\gamma(0)=x_0$ and $\gamma(1)=x$.  Since 
$(n_{i_0})_{x_0}\notin{\it l}$, $\gamma$ is $D_{\it l}$-horizontal.  
For the holomorphic isometries $h^{D_{\it l}}_{\gamma,1}:
(W_{\it l})_{x_0}\to(W_{\it l})_x$ and 
$h^{{\widehat D}_{\it l}}_{\gamma,1}
:({\widehat W}_{\it l})_{x_0}\to({\widehat W}_{\it l})_x$ as in 
Lemma 4.5, we have 
$h^{D_{\it l}}_{\gamma,1}(L^{D_{\it l}}_{x_0})=L^{D_{\it l}}_x$ 
and $h^{{\widehat D}_{\it l}}_{\gamma,1}(L^{{\widehat D}_{\it l}}_{x_0})
=L^{{\widehat D}_{\it l}}_x$.  
On the other hand, in case of $i_0\not=0$, we can show 
$h^{D_{\it l}}_{\gamma,1}=h^{{\widehat D}_{\it l}}_{\gamma,1}$ 
by imitating the discussion from Line 7 from bottom of Page 164 to Line 4 of Page 165 
in [HL2].  Also, in case of $i_0=0$, we can show 
$h^{D_{\it l}}_{\gamma,1}=h^{{\widehat D}_{\it l}}_{\gamma,1}$ 
by imitating the discussion from Line 18 of Page 165 to Line 6 of Page 166 in [HL2].  
Hence we obtain $L^{D_{\it l}}_x=L^{{\widehat D}_{\it l}}_x$.  
Therefore we obtain 
$$T_xN=\overline{(E_0)_x\oplus
\left(\mathop{\oplus}_{i\in I}(E_i)_x\right)}
=\overline{\sum\limits_{{\it l}\not={\it l}_0}T_xL^{D_{\it l}}_x}=
\overline{\sum\limits_{{\it l}\not={\it l}_0}T_xL^{\widehat D_{\it l}}_x}
=\overline{(\widehat E_0)_x\oplus
\left(\mathop{\oplus}_{i\in\widehat I}(\widehat E_i)_x\right)}=T_x\widehat N.$$
This completes the proof.  \hspace{9.45truecm}q.e.d.

\vspace{0.5truecm}

In similar to Lemma 4.2 in [HL2], we have the following fact.  

\vspace{0.5truecm}

\noindent
{\bf Lemma 4.6.3.} {\sl Let $N$ be a principal orbit of an aks-representation 
(which is a full irreducible anti-Kaehler isoparametric submanifold with $J$-diagonalizable 
shape operators).  Then each holomorphic isometry of the ambient (finite dimensional) 
anti-Kaehler space defined for $N$ in similar to the holomorphic isometry $F_{\gamma}$ 
preserves $N$ invariantly.}

\vspace{0.5truecm}

\noindent
{\it Proof.} 
Let $G/K$ be an irreducible anti-Kaehler symmetric space and $(\mathfrak g,\tau)$ the 
anti-Kaehler symmetric Lie algebra associated with $G/K$.  
Set $\mathfrak p:={\rm Ker}(\tau+{\rm id})$.  
Let $\mathfrak a_s$ be a maximal split abelian subspace of $\mathfrak p$ and 
$\mathfrak p=\mathfrak p_0+\sum\limits_{\alpha\in \triangle_+}\mathfrak p_{\alpha}$ the root 
space decomposition with respect to $\mathfrak a_s$.  
Set $\mathfrak a:=\mathfrak p_0\,(\supset\mathfrak a_s)$.  
Let $N$ be the principal orbit of the aks-representation 
$\rho:={\rm Ad}_G\vert_{\mathfrak p}:K\to GL(\mathfrak p)$ 
through a regular element $x(\in \mathfrak a)$.  Denote by $A$ the shape tensor of $N$.  
Take $v\in T^{\perp}_xN(=\mathfrak a)$ and let $\widetilde v$ be the parallel normal vector 
field of $N$ with $\widetilde v_x=v$.  Note that $\widetilde v_{\rho(k)(x)}
=\rho(k)_{\ast x}(v)$ holds for any $k\in K$.  Then we have 
$T_xN=\sum\limits_{\alpha\in\triangle_+}\mathfrak p_{\alpha}$ and 
$$A_{\widetilde v_{\rho(k)(x)}}
\vert_{\rho(k)_{\ast x}(\mathfrak p_{\alpha})}=-\frac{\alpha^{\bf c}(v)}{\alpha^{\bf c}(x)}
{\rm id}\,\,(\alpha\in\triangle_+).\leqno{(4.1)}$$
For each $\alpha\in\triangle_+$, define the section $\lambda_{\alpha}$ of 
the ${\bf C}$-dual bundle $(T^{\perp}N)^{\ast_{\bf c}}$ of $T^{\perp}N$ by 
$$(\lambda_{\alpha})_{\rho(k)(x)}:=-\frac{\alpha^{\bf c}\circ\rho(k)_{\ast x}^{-1}}
{\alpha^{\bf c}(x)}\quad(k\in K).$$
Since $\rho(k)_{\ast x}$ is the parallel translation along 
any curve $c$ in $N$ connecting $x$ and $\rho(k)(x)$ with respect to the normal connection 
of $N$, $\lambda_{\alpha}$ is a parallel section of $(T^{\perp}N)^{\ast_{\bf c}}$.  
It follows from $(4.1)$ that $\{\lambda_{\alpha}\,\vert\,
\alpha\in\triangle_+\}$ is the set of all $J$-principal curvatures of $N$.  
Let $E_{\alpha}$ be the $J$-curvature distribution for $\lambda_{\alpha}$.  
Take $\alpha_0\in\triangle_+$ and $v_0\in(\lambda_{\alpha_0})_x^{-1}(1)\setminus
\displaystyle{(\mathop{\cup}_{\alpha\in\triangle_+\,\,{\rm s.t.}\,\,
\alpha\not=\alpha_0}(\lambda_{\alpha})_x^{-1}(1))}$ and set 
$F:=\rho(K)\cdot(x+v_0)$.  It is clear that $F$ is a focal submanifold of $N$ 
whose corresponding focal distribution is equal to $E_{\alpha_0}$.  
Denote by $K_x$ (resp. $K_{x+v_0}$) the isotropy group of the $\rho(K)$-action at $x$ 
(resp. $x+v_0$) and $\mathfrak k_x$ (resp. $\mathfrak k_{x+v_0}$) the Lie algebra of 
$K_x$ (resp. $K_{x+v_0}$).  
The restriction of the $\rho(K_{x+v_0})$-action to $T^{\perp}_{x+v_0}F$ is called 
the slice representation of the $\rho(K)$-action at $x+v_0$.  It is shown that this slice 
representation coincides with the normal holonomy group action of $F$ at $x+v_0$ and that 
$\rho(K_{x+v_0})\cdot x$ is equal to $L_x^{E_{\alpha_0}}$.  
Set $\Psi(x+v_0):=\rho(K_{x+v_0})$ and $\Psi(x):=\rho(K_x)$.  
The leaf $L^{E_{\alpha_0}}_x$ is identified with the quotient manifold 
$\Psi(x+v_0)/\Psi(x)$.  
Take $X(={\rm ad}_{\mathfrak g}(\overline X))\in
{\rm Lie}\,\Psi(x+v_0)\ominus{\rm Lie}\,\Psi(x)$, where $\overline X\in
\mathfrak k_{x+v_0}$, and set $g(t):=\exp_{\Psi(x+v_0)}(tX)$ and 
$\gamma(t):=g(t)\cdot x$, where $t\in[0,1]$.  
Let $F_{\gamma}$ be the holomorphic isometry of the ambient anti-Kaehler 
space satisfying $F_{\gamma}(x)=\gamma(1)$, $(F_{\gamma})_{\ast x}
\vert_{(E_{\alpha_0})_x}
=g(1)_{\ast x}\vert_{(E_{\alpha_0})_x}$, 
$(F_{\gamma})_{\ast x}\vert_{(E_{\alpha})_x}
=h_{\gamma,1}^{E_{\alpha}}
\vert_{(E_{\alpha})_x}$ 
($\alpha\in\triangle_+\,\,{\rm s.t.}\,\,\alpha\not=\alpha_0$) and 
$(F_{\gamma})_{\ast x}\vert_{T^{\perp}_xN}=\tau^{\perp}_{\gamma}$, where 
$h_{\gamma,1}^{E_{\alpha}}$ is the holomorphic isometry defined in similar to 
$h^{D_P}_{\gamma,t}$ in the statement of Lemma 4.5 and 
$\tau^{\perp}_{\gamma}$ is the parallel translation along $\gamma$ with 
respect to the normal connection of $N$.  
Easily we can show $h_{\gamma,1}^{E_{\alpha}}
\vert_{(E_{\alpha})_x}=g(1)_{\ast x}
\vert_{(E_{\alpha})_x}$ and 
$\tau^{\perp}_{\gamma}=g(1)_{\ast x}\vert_{T_x^{\perp}N}$.  
Hence we have $(F_{\gamma})_{\ast x}=g(1)_{\ast x}$.  
Furthermore, since both $F_{\gamma}$ and $g(1)$ are affine transformations of 
the ambient anti-Kaehler space, they coincide with each other.  
Therefore, we obtain $F_{\gamma}(N)=g(1)(\rho(K)\cdot x)
=\rho(\exp_G(\overline X))(\rho(K)\cdot x)=N$.  
This completes the proof.  
\begin{flushright}q.e.d.\end{flushright}


By using Lemmas 4.6.1 and 4.6.3, we shall prove Proposition 4.6.  

\vspace{0.5truecm}

\noindent
{\it Proof of Proposition 4.6.} 
Since $M$ is a full irreducible anti-Kaehler isoparametric submanifold with 
$J$-diagonalizable shape operators and $F_{\gamma}$ is a holomorphic isometry of $V$, 
$\widehat M:=F_{\gamma}(M)$ also is a full irreducible anti-Kaehler isoparametric 
one with $J$-diagonalizable shape operators.  
Denote by $\widehat A$ the shape tensor of $\widehat M$.  
Let $\{\widehat E_0\}\cup\{\widehat E_i\,\vert\,i\in\widehat I\}$ be the set of all 
$J$-curvature distributions on $\widehat M$ and $\widehat n_i$ the $J$-curvature 
normal corresponding to $\widehat E_i$, where $\widehat E_0$ is a distribution on 
$\widehat M$ defined by 
$\displaystyle{(\widehat E_0)_x
:=\mathop{\cap}_{v\in T^{\perp}_x\widehat M}
{\rm Ker}\,\widehat A_v}$ ($x\in\widehat M$).  
Clearly we may assume that $\widehat I=I$ and $\widehat E_i=(F_{\gamma})_{\ast}(E_i)$ 
($i\in I\cup\{0\}$).  Also we have $\gamma(1)\in M\cap\widehat M$.  
Since $(F_{\gamma})_{\ast\gamma(0)}((n_i)_{\gamma(0)})
=\tau^{\perp}_{\gamma}((n_i)_{\gamma(0)})=(n_i)_{\gamma(1)}$ 
($i\in I$), we have $(\widehat n_i)_{\gamma(1)}=(n_i)_{\gamma(1)}$ ($i\in I\cup\{0\}$).  
Also, since $(F_{\gamma})_{\ast\gamma(0)}((E_i)_{\gamma(0)})
=(h^{E_i}_{\gamma,1})_{\ast\gamma(0)}((E_i)_{\gamma(0)})=(E_i)_{\gamma(1)}$ 
($i\in(I\cup\{0\})\setminus\{i_0\}$), 
we have $(\widehat E_i)_{\gamma(1)}=(E_i)_{\gamma(1)}$ ($i\in(I\cup\{0\})\setminus\{i_0\}$). 
Also, since $(F_{\gamma})_{\ast\gamma(0)}((E_{i_0})_{\gamma(0)})
=g(1)_{\ast\gamma(0)}((E_{i_0})_{\gamma(0)})=(E_{i_0})_{\gamma(1)}$, 
we have $(\widehat E_{i_0})_{\gamma(1)}=(E_{i_0})_{\gamma(1)}$.  
From these facts, we have $L^{\widehat E_i}_{\gamma(1)}=L^{E_i}_{\gamma(1)}$ 
($i\in I\cup\{0\}$) and $T_{\gamma(1)}M=T_{\gamma(1)}\widehat M$.  
Let ${\it l}_0$ be the complex affine line through $0$ and 
$(n_{i_0})_{\gamma(1)}$.  
Take any complex affine line ${\it l}$ of $T^{\perp}_{\gamma(1)}M$ 
with ${\it l}\not={\it l}_0$.  
Now we shall show that $L_{\gamma(1)}^{D_{\it l}}=L_{\gamma(1)}^{\widehat D_{\it l}}$, 
where $D_{\it l}$ (resp. $\widehat D_{\it l}$) is the distribution on $M$ 
(resp. $\widehat M$) defined as above for ${\it l}$.  
If $(n_{i_0})_{\gamma(1)}\notin{\it l}$, then 
$\gamma$ is a $D_{\it l}$-horizontal curve and hence we have 
$F_{\gamma}(L^{D_{\it l}}_{x_0})=h^{D_{\it l}}_{\gamma,1}(L^{D_{\it l}}_{x_0})
=L^{D_{\it l}}_{\gamma(1)}$ and hence $L^{\widehat D_{\it l}}_{\gamma(1)}
=L^{D_{\it l}}_{\gamma(1)}$.  
Next we consider the case of $(n_{i_0})_{\gamma(1)}\in{\it l}$.  
Then we have $0\notin{\it l}$.  
If there does not exist $i_1(\not=i_0)\in I$ with 
$(n_{i_1})_{\gamma(1)}\in{\it l}$, then 
we have $L_{\gamma(1)}^{D_{\it l}}=L_{\gamma(1)}^{E_{i_0}}
=L_{\gamma(1)}^{\widehat E_{i_0}}=L_{\gamma(1)}^{\widehat D_{\it l}}$.  
Next we consider the case where there exists $i_1(\not=i_0)\in I$ 
with $(n_{i_1})_{\gamma(1)}\in{\it l}$.  
Let $\widetilde v$ be a focal normal vector 
field of $M$ such that the corresponding focal 
distribution is equal to $D_{\it l}$.  Since $0\notin{\it l}$, it follows from 
Theorem 4.4 that $L^{D_{\it l}}_{\gamma(1)}$ is a principal orbit of the 
direct sum representation 
of aks-representations and a trivial representation.  
Since 
$(n_{i_0})_{\gamma(1)},(n_{i_1})_{\gamma(1)}\in{\it l}$ and $0\notin{\it l}$, 
$(n_{i_0})_{\gamma(1)}$ and $(n_{i_1})_{\gamma(1)}$ are ${\bf C}$-linearly 
independent.  
Assume that $L^{D_{\it l}}_{\gamma(1)}$ is reducible.  
Then the complex Coxeter group associated with $L^{D_{\it l}}_{\gamma(1)}$ is 
decomposable.  Hence $(n_{i_0})_{\gamma(1)}$ and $(n_{i_1})_{\gamma(1)}$ are orthogonal and 
there exists no $J$-curvature normal of $L^{D_{\it l}}_{\gamma(1)}$ other than them.  
Therefore, $L^{D_{\it l}}_{\gamma(1)}$ is 
congruent to the (extrinsic) product of complex spheres $L^{E_{i_0}}_{\gamma(1)}$ and 
$L^{E_{i_1}}_{\gamma(1)}$.  
Similarly $L^{D_{\it l}}_{x_0}$ is congruent to the (extrinsic) product of 
$L^{E_{i_0}}_{x_0}$ and $L^{E_{i_1}}_{x_0}$.  
Therefore we have 
$$F_{\gamma}(L^{D_{\it l}}_{x_0})=F_{\gamma}(L^{E_{i_0}}_{x_0})\times 
F_{\gamma}(L^{E_{i_1}}_{x_0})=L^{E_{i_0}}_{\gamma(1)}\times 
L^{E_{i_1}}_{\gamma(1)}=L^{D_{\it l}}_{\gamma(1)}$$
and hence $L^{\widehat D_{\it l}}_{\gamma(1)}=L^{D_{\it l}}_{\gamma(1)}$.  
Assume that $L^{D_{\it l}}_{\gamma(1)}$ is irreducible.  
Then $L^{D_{\it l}}_{\gamma(1)}$ is a principal orbit of an aks-representation.  
Then it follows from Lemma 4.6.3 that $F_{\gamma}(L^{D_{\it l}}_{x_0})=
(F_{\gamma}\vert_{(W_{\it l})_{x_0}})(L^{D_{\it l}}_{x_0})
=L^{D_{\it l}}_{x_0}$.  Hence we obtain $L^{\widehat D_{\it l}}_{\gamma(1)}
=L^{D_{\it l}}_{\gamma(1)}$.  
Thus we obtain $L^{\widehat D_{\it l}}_{\gamma(1)}=L^{D_{\it l}}_{\gamma(1)}$ 
in general.  Therefore, from Lemma 4.6.1, we obtain $M=\widehat M=F_{\gamma}(M)$, that is, 
$F_{\gamma}(M)=M$.  
\begin{flushright}q.e.d.\end{flushright}

\vspace{0.5truecm}

By using Proposition 4.6, we prove the following fact.  

\vspace{0.5truecm}

\noindent
{\bf Proposition 4.7.} {\sl For any $x\in Q(x_0)$, there exists a holomorphic 
isometry $f$ of $V$ such that $f(x_0)=x,\,f(M)=M$, $f_{\ast}(E_i)=E_i\,(i\in I)$, 
$f(Q(x_0))=Q(x_0)$ and that $f_{\ast x_0}\vert_{T^{\perp}_{x_0}M}$ coincides with the 
parallel translation along a curve in $M$ connecting $x_0$ and $x$ with respect to 
the normal connection of $M$.}

\vspace{0.5truecm}

\noindent
{\it Proof.} Take a sequence $\{x_0,x_1,\cdots,x_k(=x)\}$ of $Q(x_0)$ 
such that, for each $i\in\{0,1,\cdots,k-1\}$, $x_i$ and $x_{i+1}$ belong to a 
complex curvature sphere $S^{\bf c}_i$ of $M$.  Furthermore, for each 
$i\in\{0,1,\cdots,k-1\}$, we take the geodesic $\gamma_i:[0,1]\to S^{\bf c}_i$ with 
$\gamma_i(0)=x_i$ and $\gamma_i(1)=x_{i+1}$.  Set $f:=F_{\gamma_{k-1}}\circ\cdots\circ 
F_{\gamma_1}\circ F_{\gamma_0}$, where $F_{\gamma_i}$ ($i=0,1,\cdots,k-1$) are holomorphic 
isometries of $V$ defined in similar to the above $F_{\gamma}$.  According to 
Proposition 4.6, $f$ preserves $M$ invariantly, $f_{\ast}(E_i)=E_i\,(i\in I)$ and the 
restriction of $f_{\ast x_0}$ to $T^{\perp}_{x_0}M$ coincides with the parallel translation 
along a curve in $M$ connecting $x_0$ and $x$ with respect to the normal connection of $M$.  
Also, since $f$ preserves complex curvature spheres invariantly, it is 
shown that $f$ preserves $Q(x_0)$ invariantly.  Thus $f$ is the desired holomorphic isometry.
\begin{flushright}q.e.d.\end{flushright}

\vspace{0.5truecm}

By using Propositions 4.2 and 4.7, we shall prove Theorem A.  

\vspace{0.5truecm}

\noindent
{\it Proof of Theorem A.} Take any $\widehat x\in M$.  
Since $\overline{Q(x_0)}=M$ by Proposition 4.2, there exists a sequence 
$\{x_k\}_{k=1}^{\infty}$ in $Q(x_0)$ with 
$\lim\limits_{k\to\infty}x_k=\widehat x$.  According to Proposition 4.7, 
for each $k\in{\bf N}$, there exists a holomorphic isometry $f_k$ of $V$ with 
$f_k(x_0)=x_k,\,f_k(M)=M,\,f_k(Q(x_0))=Q(x_0)$ and 
$f_k(L^{E_i}_{x_0})=L^{E_i}_{x_k}$ ($i\in I$).  

(Step I) In this step, we shall show that, for each $i\in I$, there exists a 
subsequence $\{f_{k_j}\}_{j=1}^{\infty}$ of $\{f_k\}_{k=1}^{\infty}$ such that 
$\{f_{k_j}\vert_{L^{E_i}_{x_0}}\}_{j=1}^{\infty}$ 
pointwisely converges to a holomorphic isometry of $L^{E_i}_{x_0}$ onto 
$L^{E_i}_{\widehat x}$.  For any point $x$ of $M$, denote by $(L^{E_i}_x)_{\bf R}$ the 
compact real form through $x$ of the complex sphere $L^{E_i}_x$ satisfying 
$\langle T_x(L^{E_i}_x)_{\bf R},JT_x(L^{E_i}_x)_{\bf R}\rangle=0$, where a 
real form of $L^{E_i}_x$ means the fixed point set of an anti-holomorphic diffeomorphism 
of $L^{E_i}_x$.  Note that such a compact real form $(L^{E_i}_x)_{\bf R}$ of 
$L^{E_i}_x$ is determined uniquely (see Figure 3) and that it is isometric to 
a $m_i$-dimensional sphere, where $m_i:={\rm dim}_{\bf c}E_i$.  
Clearly we have $f_k((L^{E_i}_{x_0})_{\bf R})=(L^{E_i}_{x_k})_{\bf R}$.  
Denote by $\mathfrak F_i$ the foliation on $M$ whose leaf through $x\in M$ is equal to 
$(L^{E_i}_x)_{\bf R}$.  
Take a $\mathfrak F_i$-saturated tubular neighborhood $U$ of 
$(L^{E_i}_{\widehat x})_{\bf R}$ in $M$, where "$\mathfrak F_i$-saturatedness" of $U$ 
means that $(L^{E_i}_x)_{\bf R}\subset U$ for any $x\in U$.  Take a base $\{e_1,\cdots,e_{m_i}\}$ of 
$T_{x_0}((L^{E_i}_{x_0})_{\bf R})$ such that the norms 
$\vert\vert e_1\vert\vert,\cdots,\vert\vert e_{m_i}\vert\vert$ are sufficiently small 
and set $\bar x_a:=\exp_{x_0}(e_a)$ ($a=1,\cdots,m_i$), where $\exp_{x_0}$ is the exponential map of $(L^{E_i}_{x_0})_{\bf R}$ at $x_0$.  
Since $(L^{E_i}_x)_{\bf R}$'s ($x\in U$) are compact, 
$\mathfrak F_i$ is a Hausdorff foliation.  From this fact and the compactness of 
$(L^{E_i}_{\widehat x})_{\bf R}$, it follows that there exists a subsequence 
$\{f_{k_j}\}_{j=1}^{\infty}$ of $\{f_k\}_{k=1}^{\infty}$ such that 
$\{f_{k_j}(x_0)\}_{j=1}^{\infty}$ and $\{f_{k_j}(\bar x_a)\}_{j=1}^{\infty}$ 
($a=1,\cdots,m_i$) converge.  Set $\widehat x:=\lim\limits_{j\to\infty}f_{k_j}(x_0)$ and 
$\widehat x_a:=\lim\limits_{j\to\infty}f_{k_j}(\bar x_a)$ ($a=1,\cdots,m_i$).  
Since $\lim\limits_{j\to\infty}f_{k_j}(x_0)=\widehat x$ and 
$f_{k_j}((L^{E_i}_{x_0})_{\bf R})=(L^{E_i}_{x_{k_j}})_{\bf R}$, 
it follows from the Hausdorffness of $\mathfrak F_i$ that 
$\widehat x_a$ belongs to $(L^{E_i}_{\widehat x})_{\bf R}$ 
($a=1,\cdots,m_i$).  
Denote by $d_0,d_j$ ($j\in{\bf N}$) and $\widehat d$ the (Riemannian) distance 
functions of $(L^{E_i}_{x_0})_{\bf R},\,(L^{E_i}_{x_{k_j}})_{\bf R}$ and 
$(L^{E_i}_{\widehat x})_{\bf R}$, respectively.  Since each 
$f_{k_j}\vert_{(L^{E_i}_{x_0})_{\bf R}}$ is an isometry onto 
$(L^{E_i}_{x_{k_j}})_{\bf R}$, we have $d_j(f_{k_j}(x_0),f_{k_j}(\bar x_a))
=d_0(x_0,\bar x_a)$ and $d_j(f_{k_j}(\bar x_a),f_{k_j}(\bar x_b))=d_0(\bar x_a,\bar x_b)$, 
($a,b=1,\cdots,m_i$).  Hence we have $\widehat d(\widehat x,\widehat x_a)
=d_0(x_0,\bar x_a)$ and $\widehat d(\widehat x_a,\widehat x_b)
=d_0(\bar x_a,\bar x_b)$ ($a,b=1,\cdots,m_i$).  Therefore, since 
$(L^{E_i}_{x_0})_{\bf R}$ and $(L^{E_i}_{\widehat x})_{\bf R}$ are spheres 
isometric to each other,  there exists a unique isometry $\bar f$ 
of $(L^{E_i}_{x_0})_{\bf R}$ onto $(L^{E_i}_{\widehat x})_{\bf R}$ 
satisfying $\bar f(x_0)=\widehat x$ and $\bar f(\bar x_a)=\widehat x_a$ 
($a=1,\cdots,m_i$).  It is clear that $\bar f$ is uniquely extended to 
a holomorphic isometry of $L^{E_i}_{x_0}$ onto $L^{E_i}_{\widehat x}$.  
Denote by $f$ this holomorphic extension.  
It is easy to show that $\{f_{k_j}\vert_{(L^{E_i}_{x_0})_{\bf R}}\}_{j=1}^{\infty}$ 
pointwisely converges to $\bar f$.  Furthermore, it follows from this fact that 
$\{f_{k_j}\vert_{L^{E_i}_{x_0}}\}_{j=1}^{\infty}$ pointwisely converges to $f$.  


\vspace{0.5truecm}

\centerline{
\unitlength 0.1in
\begin{picture}( 66.9300, 28.4700)(  0.6700,-35.1000)
%
\special{pn 8}%
\special{ar 1666 1646 850 1214  5.3687718 6.2831853}%
\special{ar 1666 1646 850 1214  0.0000000 0.9056936}%
%
\special{pn 8}%
\special{ar 4488 1624 850 1214  2.2357680 4.0561528}%
%
\special{pn 13}%
\special{ar 3082 1638 558 126  6.2831853 6.2831853}%
\special{ar 3082 1638 558 126  0.0000000 3.1415927}%
%
\special{pn 8}%
\special{ar 3082 1638 558 124  3.1415927 3.1767833}%
\special{ar 3082 1638 558 124  3.2823551 3.3175457}%
\special{ar 3082 1638 558 124  3.4231176 3.4583082}%
\special{ar 3082 1638 558 124  3.5638800 3.5990707}%
\special{ar 3082 1638 558 124  3.7046425 3.7398331}%
\special{ar 3082 1638 558 124  3.8454050 3.8805956}%
\special{ar 3082 1638 558 124  3.9861674 4.0213580}%
\special{ar 3082 1638 558 124  4.1269299 4.1621205}%
\special{ar 3082 1638 558 124  4.2676924 4.3028830}%
\special{ar 3082 1638 558 124  4.4084548 4.4436454}%
\special{ar 3082 1638 558 124  4.5492173 4.5844079}%
\special{ar 3082 1638 558 124  4.6899798 4.7251704}%
\special{ar 3082 1638 558 124  4.8307422 4.8659328}%
\special{ar 3082 1638 558 124  4.9715047 5.0066953}%
\special{ar 3082 1638 558 124  5.1122671 5.1474578}%
\special{ar 3082 1638 558 124  5.2530296 5.2882202}%
\special{ar 3082 1638 558 124  5.3937921 5.4289827}%
\special{ar 3082 1638 558 124  5.5345545 5.5697451}%
\special{ar 3082 1638 558 124  5.6753170 5.7105076}%
\special{ar 3082 1638 558 124  5.8160795 5.8512701}%
\special{ar 3082 1638 558 124  5.9568419 5.9920325}%
\special{ar 3082 1638 558 124  6.0976044 6.1327950}%
\special{ar 3082 1638 558 124  6.2383668 6.2735575}%
%
\special{pn 20}%
\special{sh 1}%
\special{ar 3082 1646 10 10 0  6.28318530717959E+0000}%
\special{sh 1}%
\special{ar 3082 1646 10 10 0  6.28318530717959E+0000}%
%
\special{pn 8}%
\special{pa 2898 1756}%
\special{pa 3280 1528}%
\special{dt 0.045}%
%
\special{pn 20}%
\special{sh 1}%
\special{ar 2906 1764 10 10 0  6.28318530717959E+0000}%
\special{sh 1}%
\special{ar 2906 1764 10 10 0  6.28318530717959E+0000}%
%
\special{pn 20}%
\special{sh 1}%
\special{ar 3272 1528 10 10 0  6.28318530717959E+0000}%
\special{sh 1}%
\special{ar 3272 1528 10 10 0  6.28318530717959E+0000}%
%
\special{pn 8}%
\special{ar 2632 1764 268 1200  5.2181430 6.2831853}%
\special{ar 2632 1764 268 1200  0.0000000 0.9505468}%
%
\special{pn 13}%
\special{pa 3680 1950}%
\special{pa 3650 1960}%
\special{pa 3618 1968}%
\special{pa 3586 1970}%
\special{pa 3554 1972}%
\special{pa 3522 1972}%
\special{pa 3490 1970}%
\special{pa 3460 1968}%
\special{pa 3428 1964}%
\special{pa 3396 1958}%
\special{pa 3364 1954}%
\special{pa 3332 1948}%
\special{pa 3302 1940}%
\special{pa 3270 1932}%
\special{pa 3240 1922}%
\special{pa 3210 1914}%
\special{pa 3180 1904}%
\special{pa 3150 1892}%
\special{pa 3120 1880}%
\special{pa 3090 1868}%
\special{pa 3062 1854}%
\special{pa 3032 1842}%
\special{pa 3004 1826}%
\special{pa 2976 1810}%
\special{pa 2950 1792}%
\special{pa 2924 1772}%
\special{pa 2898 1754}%
\special{pa 2898 1754}%
\special{sp}%
%
\special{pn 8}%
\special{ar 2656 1888 800 574  4.4960064 4.5134864}%
\special{ar 2656 1888 800 574  4.5659263 4.5834063}%
\special{ar 2656 1888 800 574  4.6358462 4.6533262}%
\special{ar 2656 1888 800 574  4.7057661 4.7232460}%
\special{ar 2656 1888 800 574  4.7756860 4.7931659}%
\special{ar 2656 1888 800 574  4.8456058 4.8630858}%
\special{ar 2656 1888 800 574  4.9155257 4.9330057}%
\special{ar 2656 1888 800 574  4.9854456 5.0029256}%
\special{ar 2656 1888 800 574  5.0553655 5.0728455}%
\special{ar 2656 1888 800 574  5.1252854 5.1427653}%
\special{ar 2656 1888 800 574  5.1952053 5.2126852}%
\special{ar 2656 1888 800 574  5.2651251 5.2826051}%
\special{ar 2656 1888 800 574  5.3350450 5.3525250}%
\special{ar 2656 1888 800 574  5.4049649 5.4224449}%
\special{ar 2656 1888 800 574  5.4748848 5.4923648}%
\special{ar 2656 1888 800 574  5.5448047 5.5622846}%
%
\special{pn 8}%
\special{pa 3260 1516}%
\special{pa 3288 1532}%
\special{pa 3316 1548}%
\special{pa 3344 1564}%
\special{pa 3370 1582}%
\special{pa 3396 1600}%
\special{pa 3422 1618}%
\special{pa 3446 1638}%
\special{pa 3470 1660}%
\special{pa 3494 1682}%
\special{pa 3516 1704}%
\special{pa 3538 1728}%
\special{pa 3558 1752}%
\special{pa 3580 1776}%
\special{pa 3598 1802}%
\special{pa 3616 1830}%
\special{pa 3632 1856}%
\special{pa 3646 1886}%
\special{pa 3660 1914}%
\special{pa 3674 1944}%
\special{pa 3680 1960}%
\special{sp -0.045}%
%
\special{pn 13}%
\special{pa 2496 1336}%
\special{pa 2498 1368}%
\special{pa 2508 1398}%
\special{pa 2520 1428}%
\special{pa 2536 1456}%
\special{pa 2554 1482}%
\special{pa 2572 1508}%
\special{pa 2592 1532}%
\special{pa 2614 1556}%
\special{pa 2636 1580}%
\special{pa 2660 1602}%
\special{pa 2684 1622}%
\special{pa 2708 1642}%
\special{pa 2734 1662}%
\special{pa 2760 1680}%
\special{pa 2788 1698}%
\special{pa 2814 1714}%
\special{pa 2842 1730}%
\special{pa 2872 1744}%
\special{pa 2900 1756}%
\special{pa 2920 1764}%
\special{sp}%
%
\special{pn 8}%
\special{pa 1998 1662}%
\special{pa 2582 1522}%
\special{fp}%
\special{sh 1}%
\special{pa 2582 1522}%
\special{pa 2512 1518}%
\special{pa 2530 1534}%
\special{pa 2522 1556}%
\special{pa 2582 1522}%
\special{fp}%
%
\special{pn 8}%
\special{pa 2008 2006}%
\special{pa 2632 1712}%
\special{fp}%
\special{sh 1}%
\special{pa 2632 1712}%
\special{pa 2562 1722}%
\special{pa 2584 1736}%
\special{pa 2580 1760}%
\special{pa 2632 1712}%
\special{fp}%
\put(19.5700,-19.6900){\makebox(0,0)[rt]{$(L^{E_i}_x)_{\bf R}$}}%
\put(19.8100,-16.1700){\makebox(0,0)[rt]{$L'$}}%
\put(20.0000,-29.9000){\makebox(0,0)[lt]{$L'$ is a compact real form of $L^{E_i}_x$ }}%
\put(20.1000,-31.9000){\makebox(0,0)[lt]{but $\langle T_xL',JT_xL'\rangle\not=0$}}%
\put(42.5000,-7.5000){\makebox(0,0)[lt]{$L^{E_i}_x$}}%
%
\special{pn 8}%
\special{pa 4216 854}%
\special{pa 3532 1294}%
\special{fp}%
\special{sh 1}%
\special{pa 3532 1294}%
\special{pa 3598 1274}%
\special{pa 3576 1264}%
\special{pa 3576 1240}%
\special{pa 3532 1294}%
\special{fp}%
%
\special{pn 8}%
\special{pa 1980 2450}%
\special{pa 2900 1770}%
\special{fp}%
\special{sh 1}%
\special{pa 2900 1770}%
\special{pa 2836 1794}%
\special{pa 2858 1802}%
\special{pa 2858 1826}%
\special{pa 2900 1770}%
\special{fp}%
\put(19.3000,-24.2000){\makebox(0,0)[rt]{$x$}}%
%
\special{pn 8}%
\special{pa 5100 1430}%
\special{pa 6400 1430}%
\special{pa 6400 2580}%
\special{pa 5100 2580}%
\special{pa 5100 1430}%
\special{fp}%
%
\special{pn 20}%
\special{sh 1}%
\special{ar 5750 2010 10 10 0  6.28318530717959E+0000}%
\special{sh 1}%
\special{ar 5750 2010 10 10 0  6.28318530717959E+0000}%
\put(55.4000,-27.6000){\makebox(0,0)[lt]{$T_xL^{E_i}_x$}}%
%
\special{pn 8}%
\special{pa 5100 2010}%
\special{pa 6400 2010}%
\special{fp}%
%
\special{pn 8}%
\special{pa 5740 1430}%
\special{pa 5740 2580}%
\special{fp}%
%
\special{pn 8}%
\special{pa 5100 1810}%
\special{pa 6400 2200}%
\special{fp}%
%
\special{pn 8}%
\special{pa 5550 2580}%
\special{pa 5940 1430}%
\special{fp}%
%
\special{pn 8}%
\special{pa 6570 1720}%
\special{pa 6240 2010}%
\special{dt 0.045}%
\special{sh 1}%
\special{pa 6240 2010}%
\special{pa 6304 1982}%
\special{pa 6280 1976}%
\special{pa 6278 1952}%
\special{pa 6240 2010}%
\special{fp}%
%
\special{pn 8}%
\special{pa 5740 1190}%
\special{pa 5910 1520}%
\special{dt 0.045}%
\special{sh 1}%
\special{pa 5910 1520}%
\special{pa 5898 1452}%
\special{pa 5886 1474}%
\special{pa 5862 1470}%
\special{pa 5910 1520}%
\special{fp}%
%
\special{pn 8}%
\special{pa 6670 1990}%
\special{pa 6310 2160}%
\special{dt 0.045}%
\special{sh 1}%
\special{pa 6310 2160}%
\special{pa 6380 2150}%
\special{pa 6358 2138}%
\special{pa 6362 2114}%
\special{pa 6310 2160}%
\special{fp}%
\put(55.1000,-12.5000){\makebox(0,0)[rb]{$JT_x(L^{E_i}_x)_{\bf R}$}}%
\put(65.9000,-16.3000){\makebox(0,0)[lt]{$T_x(L^{E_i}_x)_{\bf R}$}}%
\put(66.9000,-19.2000){\makebox(0,0)[lt]{$T_xL'$}}%
\put(55.7000,-11.6000){\makebox(0,0)[lb]{$JT_xL'$}}%
\put(50.1000,-30.8000){\makebox(0,0)[lt]{$\langle T_x(L^{E_i}_x)_{\bf R},JT_x(L^{E_i}_x)_{\bf R}\rangle=0$}}%
\put(50.2000,-33.0000){\makebox(0,0)[lt]{$\langle T_xL',JT_xL'\rangle\not=0$}}%
%
\special{pn 8}%
\special{pa 5100 2580}%
\special{pa 6400 1430}%
\special{da 0.070}%
%
\special{pn 8}%
\special{pa 5100 1430}%
\special{pa 6400 2580}%
\special{da 0.070}%
%
\special{pn 8}%
\special{pa 4770 1920}%
\special{pa 5220 1540}%
\special{da 0.070}%
\special{sh 1}%
\special{pa 5220 1540}%
\special{pa 5156 1568}%
\special{pa 5180 1574}%
\special{pa 5182 1598}%
\special{pa 5220 1540}%
\special{fp}%
%
\special{pn 8}%
\special{pa 4770 1920}%
\special{pa 5380 2320}%
\special{da 0.070}%
\special{sh 1}%
\special{pa 5380 2320}%
\special{pa 5336 2268}%
\special{pa 5336 2292}%
\special{pa 5314 2300}%
\special{pa 5380 2320}%
\special{fp}%
\put(47.3000,-18.7000){\makebox(0,0)[rt]{{\small null cone}}}%
%
\special{pn 8}%
\special{pa 5510 1430}%
\special{pa 5980 2580}%
\special{fp}%
%
\special{pn 8}%
\special{pa 5430 1260}%
\special{pa 5740 1540}%
\special{dt 0.045}%
\special{sh 1}%
\special{pa 5740 1540}%
\special{pa 5704 1480}%
\special{pa 5700 1504}%
\special{pa 5678 1510}%
\special{pa 5740 1540}%
\special{fp}%
%
\special{pn 8}%
\special{pa 6700 2340}%
\special{pa 5940 2450}%
\special{dt 0.045}%
\special{sh 1}%
\special{pa 5940 2450}%
\special{pa 6010 2460}%
\special{pa 5994 2442}%
\special{pa 6004 2422}%
\special{pa 5940 2450}%
\special{fp}%
\put(67.6000,-22.7000){\makebox(0,0)[lt]{$V'$}}%
\put(50.2000,-35.1000){\makebox(0,0)[lt]{$\langle T_xL',V'\rangle=0$}}%
\end{picture}%
\hspace{4.7truecm}}

\vspace{1.2truecm}

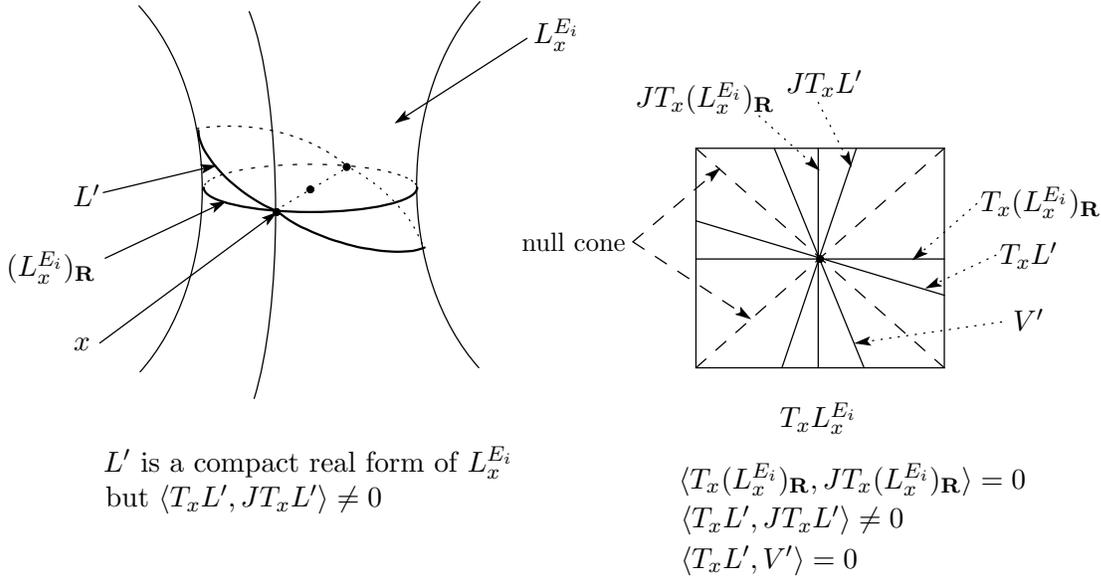
\centerline{{\bf Figure 3.}}

\vspace{0.5truecm}

(Step II) Next we shall show that, for each fixed $y\in Q(x_0)$, there exists 
a subsequence $\{f_{k_j}\}_{j=1}^{\infty}$ of $\{f_k\}_{k=1}^{\infty}$ such 
that $\{f_{k_j}(y)\}_{j=1}^{\infty}$ converges.  
There exists a sequence 
$\{\bar x_0(=x_0),\bar x_1,\cdots,\bar x_m(=y)\}$ in $Q(x_0)$ such that, 
for each $j\in\{1,\cdots,m\}$, $\bar x_j$ is contained in a complex curvature 
sphere through $\bar x_{j-1}$ (which we denote by $L^{E_{i(j)}}_{\bar x_{j-1}}$).  
For simplicity, we shall consider the 
case of $m=2$.  From the fact in Step I, 
there exists a subsequence $\{f_{k^1_j}\}_{j=1}^{\infty}$ of 
$\{f_k\}_{k=1}^{\infty}$ 
such that $\{f_{k^1_j}\vert_{L^{E_{i(1)}}_{x_0}}\}_{j=1}^{\infty}$ pointwisely 
converges to a holomorphic isometry $f^1$ of $L^{E_{i(1)}}_{x_0}$ onto 
$L^{E_{i(1)}}_{\widehat x}$.  Furthermore, by noticing 
$\lim\limits_{j\to\infty}f_{k^1_j}(\bar x_1)=f^1(\bar x_1)$ and imitating 
the discussion in Step I, we can show that there exists a subsequence 
$\{f_{k^2_j}\}_{j=1}^{\infty}$ of $\{f_{k^1_j}\}_{j=1}^{\infty}$ such that 
$\{f_{k^2_j}\vert_{L^{E_{i(2)}}_{\bar x_1}}\}_{j=1}^{\infty}$ pointwisely 
converges to a holomorphic isometry $f^2$ of $L^{E_{i(2)}}_{\bar x_1}$ onto 
$L^{E_{i(2)}}_{f^1(\bar x_1)}$.  
Since $y=\bar x_2\in L^{E_{i(2)}}_{\bar x_1}$, we have 
$\lim\limits_{j\to\infty}f_{k^2_j}(y)=f^2(y)$.  
Thus $\{f_{k^2_j}\}_{j=1}^{\infty}$ is the desired subsequence of 
$\{f_{k_j}\}_{j=1}^{\infty}$.  

(Step III) Let $W$ be the complex affine span of $M$.  
Next we shall show that there exists a subsequence 
$\{f_{k_j}\}_{j=1}^{\infty}$ of $\{f_k\}_{k=1}^{\infty}$ 
such that $\{f_{k_j}\vert_W\}_{j=1}^{\infty}$ pointwisely converges to some 
holomorphic isometry of $W$.  
Take a countable subset $B:=\{w_j\,\vert\,j\in{\bf N}\}$ of $Q(x_0)$ with 
$\overline B=\overline{Q(x_0)}(=M)$.  According to the fact in Step II, there exists 
a subsequence $\{f_{k^1_j}\}_{j=1}^{\infty}$ 
of $\{f_k\}_{k=1}^{\infty}$ such that $\{f_{k^1_j}(w_1)\}_{j=1}^{\infty}$ converges.  
Again, according to the fact in Step II, there exists a 
subsequence $\{f_{k^2_j}\}_{j=1}^{\infty}$ of $\{f_{k^1_j}\}_{j=1}^{\infty}$ 
such that $\{f_{k^2_j}(w_2)\}_{j=1}^{\infty}$ converges.  
In the sequel, we take subsequences $\{f_{k^{\it l}_j}\}_{j=1}^{\infty}$ 
(${\it l}=3,4,5,\cdots$) inductively.  It is clear that 
$\{f_{k^j_j}(w_{\it l})\}_{j=1}^{\infty}$ converges for each ${\it l}\in
{\bf N}$, that is, $\{f_{k^j_j}\vert_B\}_{j=1}^{\infty}$ pointwisely converges 
to some map $f$ of $B$ into $M$.  Since each $f_{k^j_j}$ is a holomorphic isometry and hence 
$f_{k^j_j}\vert_W:W\to W$ is an affine transformation, $f$ extends to an affine 
transformation of $W$.  Denote by $\widetilde f$ this extension.  It is clear that 
$\{f_{k^j_j}\vert_W\}_{j=1}^{\infty}$ pointwisely converges to $\widetilde f$ and that 
$\widetilde f$ is a holomorphic isometry of $W$.  

(Step IV) Denote by $H$ the group generated by all holomorphic isometries of $V$ preserving 
$M$ invariantly.  Let $\widetilde f$ be as in Step III.  
It is clear that $\widetilde f$ extends to a holomorphic isometry of $V$.  
Denote by $\widehat f$ this extension.  It is clear that $\widehat f(M)=M$ and 
$\widehat f(x_0)=\lim\limits_{j\to\infty}f_{k^j_j}(x_0)
=\lim\limits_{j\to\infty}x_{k^j_j}=\widehat x$.  
Hence we have $\widehat x\in H\cdot x_0$.  From the arbitrariness of 
$\widehat x$, we obtain $M\subset H\cdot x_0$.  
On the ther hand, it is clear that $H\cdot x_0\subset M$.  
Therefore we obtain $H\cdot x_0=M$.  
\hspace{10.7truecm}q.e.d.

\vspace{1truecm}

\centerline{{\bf References}}

\vspace{0.5truecm}

\small{
\noindent
[Be] M. Berger, Les espaces sym$\acute e$triques non compacts, 
Ann. Sci. $\acute E$c. Norm. Sup$\acute e$r. III. S$\acute e$r. {\bf 74} 

(1959) 85-177.

\noindent
[Br] M. Br$\ddot u$ck, Equifocal famlies in symmetric spaces of compact type, J.reine angew. 
Math. 

{\bf 515} (1999), 73-95.

\noindent
[BCO] J. Berndt, S. Console and C. Olmos, Submanifolds and holonomy, Research 
Notes in 

Mathematics 434, CHAPMAN $\&$ HALL/CRC Press, Boca Raton, London, New York 

Washington, 2003.

\noindent
[BH] R. A. Blumenthal and J. J. Hebda, Complementary distributions which 
preserves the leaf 

geometry and applications to totally geodesic foliations, Quat. J. 
Math. {\bf 35} (1984) 383-392.









\noindent
[Ca] E. Cartan, Familles de surfaces isoparam$\acute e$triques dans les 
espaces $\acute a$ courbure constante, 

Ann. Mat. Pura Appl. {\bf 17} (1938), 177--191.

\noindent
[CP] M. Cahen and M. Parker, Pseudo-riemannian symmetric spaces, Memoirs 
of the Amer. Math. 

Soc. {\bf 24} No. 229 (1980).  

\noindent
[Ch] U. Christ, 
Homogeneity of equifocal submanifolds, J. Differential Geom. {\bf 62} (2002), 
1--15.

\noindent
[Cox] H. S. M. Coxeter, Discrete groups generated by reflections, 
Ann. of Math. (2) {\bf 35} (1934),

588--621.

\noindent
[E] H. Ewert, 
Equifocal submanifolds in Riemannian symmetric spaces, Doctoral thesis.

\noindent
[G1] L. Geatti, 
Invariant domains in the complexfication of a noncompact Riemannian 
symmetric 

space, J. Algebra {\bf 251} (2002), 619--685.

\noindent
[G2] L. Geatti, 
Complex extensions of semisimple symmetric spaces, manuscripta math. {\bf 120} 

(2006) 1-25.

\noindent
[GG] L. Geatti and C. Gorodski, 
Polar orthogonal representations of real reductive algebraic 

groups, J. Algebra {\bf 320} (2008) 3036-3061.



\noindent
[Ha] J. Hahn, Isotropy representations of semisimple symmetric spaces 
and homogeneous hyper-

surfaces, J. Math. Soc. Japan {\bf 40} (1988) 271-288.  

\noindent
[HL1] E. Heintze and X. Liu, 
A splitting theorem for isoparametric submanifolds in Hilbert space, 

J. Differential Geom. {\bf 45} (1997), 319--335.

\noindent
[HL2] E. Heintze and X. Liu, 
Homogeneity of infinite dimensional isoparametric submanifolds, 

Ann. of Math. {\bf 149} (1999), 149-181.

\noindent
[HLO] E. Heintze, X. Liu and C. Olmos, 
Isoparametric submanifolds and a Chevalley type restric-

ction theorem, Integrable systems, geometry, and topology, 151-190, 
AMS/IP Stud. Adv. 

Math. 36, Amer. Math. Soc., Providence, RI, 2006.

\noindent
[HOT] E. Heintze, C. Olmos and G. Thorbergsson, 
Submanifolds with constant prinicipal curva-

tures and normal holonomy groups, Int. J. Math. {\bf 2} (1991), 167-175.

\noindent
[HPTT] E. Heintze, R. S. Palais, C. L. Terng and G. Thorbergsson, 
Hyperpolar actions on sym-

metric spaces, Geometry, topology and physics, 214--245 Conf. Proc. Lecture 
Notes Geom. 

Topology {\bf 4}, Internat. Press, Cambridge, Mass., 1995.

\noindent
[He] S. Helgason, 
Differential geometry, Lie groups and symmetric spaces, Pure Appl. Math. 80, 

Academic Press, New York, 1978.

\noindent
[Hu] M. C. Hughes, 
Complex reflection groups, Comm. Algebra {\bf 18} (1990), 3999--4029.

\noindent
[KN] S. Kobayashi and K. Nomizu, 
Foundations of differential geometry, Interscience Tracts in 

Pure and Applied Mathematics 15, Vol. II, New York, 1969.



\noindent
[Koi1] N. Koike, 
Submanifold geometries in a symmetric space of non-compact 
type and a pseudo-

Hilbert space, Kyushu J. Math. {\bf 58} (2004), 167--202.

\noindent
[Koi2] N. Koike, 
Complex equifocal submanifolds and infinite dimensional 
anti-Keahlerian isopa-

rametric submanifolds, Tokyo J. Math. {\bf 28} (2005), 
201--247.



\noindent
[Koi3] N. Koike, 
A splitting theorem for proper complex equifocal submanifolds, Tohoku Math. 

J. {\bf 58} (2006) 393-417.



\noindent
[Koi4] N. Koike, 
The homogeneous slice theorem for the complete complexification of a 
proper 

complex equifocal submanifold, Tokyo J. Math. {\bf 33} (2010), 1-30.





\noindent
[Koi5] N. Koike, Hermann type actions on a pseudo-Riemannian symmetric 
space, Tsukuba J. 

Math. {\bf 34} (2010), 137-172.







\noindent
[O] C. Olmos, Isoparametric submanifolds and their homogeneous structures, 
J. Differential 

Geom. {\bf 38} (1993), 225-234.

\noindent
[OW] C. Olmos and A. Will, Normal holonomy in Lorentzian space and 
submanifold geometry, 

Indiana Univ. Math. J. {\bf 50} (2001), 1777-1788.

\noindent
[O'N] B. O'Neill, 
Semi-Riemannian Geometry, with applications to relativity, 
Pure Appl. Math. 

103, Academic Press, New York, 1983.

\noindent
[OS] T. Ohshima and J. Sekiguchi, The restricted root system of 
a semisimple symmetric pair, 

Group representations and systems of 
differential equations (Tokyo, 1982), 433--497, Adv. 

Stud. Pure Math. {\bf 4}, North-Holland, Amsterdam, 1984.

\noindent
[Pa] R. S. Palais, 
Morse theory on Hilbert manifolds, Topology {\bf 2} (1963), 299--340.



\noindent
[PT2] R. S. Palais and C. L. Terng, Critical point theory and submanifold 
geometry, Lecture Notes 

in Math. {\bf 1353}, Springer-Verlag, Berlin, 1988.

\noindent
[Pe] A. Z. Petrov, Einstein spaces, Pergamon Press, 1969.

\noindent
[R] W. Rossmann, The structures of semisimple symmetric spaces, 
Canad. J. Math. {\bf 31} (1979), 

157-180.

\noindent
[Si] J. Simons, On the transitivity of holonomy systems, Ann. of 
Math. {\bf 76} (1962), 213--234.

\noindent
[Sz1] R. Sz$\ddot{{{\rm o}}}$ke, Involutive structures on the 
tangent bundle of symmetric spaces, Math. Ann. {\bf 319} 

(2001), 319--348.

\noindent
[Sz2] R. Sz$\ddot{{{\rm o}}}$ke, Canonical complex structures associated to 
connections and complexifications of 

Lie groups, Math. Ann. {\bf 329} (2004), 553--591.

\noindent
[Te1] C. L. Terng, 
Isoparametric submanifolds and their Coxeter groups, 
J. Differential Geom. {\bf 21} 

(1985), 79--107.

\noindent
[Te2] C. L. Terng, 
Proper Fredholm submanifolds of Hilbert space, 
J. Differential Geom. {\bf 29} 

(1989), 9--47.

\noindent
[Te3] C. L. Terng, 
Polar actions on Hilbert space, J. Geom. Anal. {\bf 5} (1995), 
129--150.

\noindent
[TT] C. L. Terng and G. Thorbergsson, 
Submanifold geometry in symmetric spaces, J. Differential 

Geom. {\bf 42} (1995), 665--718.

\noindent
[Th] G. Thorbergsson, Isoparametric 
foliations and their buildings, Ann of Math. {\bf 133} (1991), 

429--446.

\noindent
[W1] H. Wu, Holonomy groups of indefinite metrics, Pacific J. Math. {\bf 20} 
(1967), 351--392.

\noindent
[W2] B. Wu, Isoparametric submanifolds of hyperbolic spaces, 
Trans. Amer. Math. Soc. {\bf 331} 

(1992), 609--626.}



\end{document}